\documentclass{article}
\usepackage{graphicx} 

\title{Künneth and extension theorems for Morita invariants of non-commutative schemes}
\author{Alexandre Nicolle}
\date{December 2024}

\usepackage[utf8]{inputenc}
\usepackage[english]{babel}
\usepackage[T1]{fontenc}
\usepackage{amsmath}
\usepackage{amsfonts}
\usepackage{amssymb}
\usepackage{multicol}
\usepackage{array}
\usepackage{mathrsfs}
\usepackage{amsthm}
\usepackage[all]{xy} 
\usepackage{marvosym}
\usepackage{tikz-cd}

\tikzset{%
    symbol/.style={%
        draw=none,
        every to/.append style={%
            edge node={node [sloped, allow upside down, auto=false]{$#1$}}}
    }
}
\usepackage{calrsfs}

\usepackage[colorlinks=true, linkcolor=violet, citecolor=orange, anchorcolor=orange,]{hyperref}

\newtheorem{theorem}{\textsc{Theorem}}[section]
\newtheorem{lemma}{\textsc{Lemma}}[section]
\newtheorem{proposition}{\textsc{Proposition}}[section]
\newtheorem{corollary}{\textsc{Corollary}}[section]
\newtheorem{remarque}{\textsc{Remark}}[section]
\newtheorem{definition}{\textsc{Definition}}[section]

\newtheorem*{keyword}{\textsc{Keywords}}
\newtheorem{example}{\textsc{Example}}[section]

\newcommand{\eq}{\Leftrightarrow}

\usepackage{lmodern, mathtools}
\usepackage{old-arrows} 
\DeclareMathOperator\Hom{Hom}
\DeclareMathOperator\coker{coker}
\DeclareMathOperator\RHom{RHom}
\DeclareMathOperator\Ext{Ext}

\DeclareMathOperator\Hoc{HH}

\DeclareMathOperator\HC{HC}
\DeclareMathOperator\Coh{Coh}
\DeclareMathOperator\Fun{Fun}
\DeclareMathOperator\Ho{Ho}
\DeclareMathOperator\dg{{dg}}
\DeclareMathOperator\Hh{H}
\DeclareMathOperator\Ob{Ob}

\DeclareMathOperator\dgcat{dg-Cat}
\DeclareMathOperator\ac{Ac}

\DeclareMathOperator\Hqe{Ho_Q(dg-Cat_R)}
\DeclareMathOperator\op{{op}}
\DeclareMathOperator\aA{\mathcal{A}}
\DeclareMathOperator\bB{\mathcal{B}}
\DeclareMathOperator\cC{\mathcal{C}}
\DeclareMathOperator\dD{\mathcal{D}}
\DeclareMathOperator\pP{\mathcal{P}}
\DeclareMathOperator\qQ{\mathcal{Q}}
\DeclareMathOperator\sS{\mathcal{S}}
\DeclareMathOperator\tT{\mathcal{T}}
\DeclareMathOperator\uU{\mathcal{U}}

\DeclareMathOperator\vV{\mathcal{V}}

\DeclareMathOperator\xX{\mathcal{X}}
\DeclareMathOperator\yY{\mathcal{Y}}
\DeclareMathOperator\zZ{\mathcal{Z}}

\DeclareMathOperator\spec{{Spec}}

\DeclareMathOperator\id{{id}}

\DeclareMathOperator\set{{Set}}
\DeclareMathOperator\qcoh{{Qcoh}}
\DeclareMathOperator\Mod{{Mod}}

\DeclareMathOperator\qfun{{QFun}}
\DeclareMathOperator\rep{{Rep}}

\DeclareMathOperator\bif{{Int}}
\DeclareMathOperator\perf{{Perf}}
\DeclareMathOperator\hproj{h-Proj}
\DeclareMathOperator\hflat{h-Flat}
\DeclareMathOperator\Sf{SF}
\DeclareMathOperator\ra{\rightarrow}

\DeclareMathOperator\HP{HP}

\usepackage[left=3.0cm,right=3.0cm,top=2.2cm,bottom=2.5cm]{geometry}

\numberwithin{equation}{section}

\begin{document}

\maketitle

\begin{center}
    \emph{En hommage\footnote{Respectively, for his fatherhood regarding the philosophy of general position ;  and for his love of systematizing humble things until they magnify.} à René Thom et Jean Dieudonné.}\\
    
\end{center}

\begin{abstract}
    In this article, we apply the derived Morita theory of dg-categories to show how to extend the domain of validity of many identities relating Morita invariants from associative dg-algebras toward non-commutative scheme. Doing so, we obtain that the dg-category of associative algebras can be used to test the exactness of any sequence and the commutativity of any diagram involving Morita invariants depending multilinearly on their arguments, under a mild condition of stability by cofibrant replacement. This gives a simple picture of how the underlying algebra of a non-commutative scheme captures its Morita invariant properties. As an application, we use a  Künneth formula on non-commutative schemes to factorise the Hochschild cohomology of a product of quasi-phantom categories such as built by Orlov and Gorshinsky in \cite{Gorchinskiy_2013}. This is an occasion to explore how the derived Morita theory interacts with the Čech dg-enhancement of a projective scheme and to shed an unifying light upon the landscape of dg-modules over dg-categories.\\ 
\end{abstract}

\begin{keyword}Noncommutative geometry, derived algebraic geometry, dg-categories, homotopy theory, model categories, noncommutative schemes, Morita invariants, Künneth formula, phantom categories.
\end{keyword}

\tableofcontents

\paragraph{Some notations}
In all the article, "dg" will mean differential(ly) graded. We preferentially use the typography $\aA,\bB,\cdots$ for dg-categories (Definition \ref{dgcatdef}), $\mathfrak{A},\mathfrak{B},\cdots $ for dg-algebras, $x,y,\cdots $ for objects in dg-categories, $f,g$ for dg-functors. $(-)^\omega$ shall denote subcategories of compact objects defined according to the context, $\Gamma$ shall denote a generic weak Morita invariant (Definition \ref{defweakmo}), $X,Y$ for (commutative) schemes, $F,G$ for coherent sheaves over them, $R$ for a commutative ring, $\xX,\yY$ for noncommutative schemes (Definition \ref{noncomschdef}).

\section*{Introduction}

\paragraph{Context}In derived algebraic geometry, non-commutative schemes -- i.e dg-categories of perfect complexes of modules over a cohomologically bounded associative dg-algebra -- form an important class of objects. It's importance has been pointed out twenty years ago by the school of Kontsevitch \cite{Kontsevitch} in search for a good theory of "sheaves of categories", unifying certain phenomena common to derived symplectic and algebraic geometry, with application to topological field theories. In particular, derived noncommutative schemes are a suitable choice of base spaces over which to consider complexes of sheaves of modules, with the hope to characterize the geometry of the underlying space through some reconstruction theorem. The theory and main properties of non-commutative schemes have been thoroughly studied and carved-out, notably by Orlov \cite{Orlov_2016,ORLOV2020107096}, but many ground structures and ideas appears earlier under a more implicit and application minded form, to treat of associative algebras or projective schemes through dg-categorical methods : One the algebraic side of the force we shall cite the foundational work of Keller \cite{Keller1994, keller2006differentialgradedcategories, keller:hal-03280332} focusing on the perspectives for non commutative algebras themselves; and on the geometric side one of the main inspiration may be found in Kuznetsov's work \cite{kuznetsov2009, Kuznetsov_2013}, who has been exploring the use of this formalism to compute  invariants of non-commutative nature over projective schemes with great success. Derived noncommutative algebraic geometry offers a unifying perspective to define and study the cohomology of fibered algebraic spaces i.e geometrical object such as algebraic bundles over schemes, that can be associated a compactly generated triangulated structure. In particular the formalism of derived non-commutative schemes resolve many technicities arising in the previous attempts for an underived noncommutative geometry \cite{Orlov_2003,Rosenberg98}, which have been somewhat hindered by the complexity of a direct approach to a generalised categorical Grothendieck construction : In fact, if you want to classify objects in a category of spaces (in a broad sense), then you basically need to describe the moduli stack classifying those objects, and for reasons whose explanation exceeds the frame of our paper this is much more naturally done up to homotopy\footnote{In one word : The problem which makes homotopy unavoidable in algebraic geometry is that one must take into account the excess of intersection multiplicity arising while working with objects such a schemes. This can only be done by putting the spaces we want to classify in general position in order to resorb some local bad phenomenons. But this is even more generally true for the following reason, which is a broader reformulation of the previous informal argument about "bad local phenomenons" :  Studying the homotopical perturbations of a space, is exactly studying how the points associated to that space in the classifying stack glues to its neighbourhood, this means that even if you do not quotient by deformations, i.e work up to homotopy, you still have to care about homotopy, but you will also be forced to manage all pathological situations by hand, which is a huge problem, because these pathologies will cumulate with local pathological behaviours inherent to stacky points, bon courage!}. \\

\paragraph{Content}As evoked in the previous paragraph, it is a necessity to work in the derived frame while dealing with classification problems in geometry : this leads to consider replacements for the notion of homotopy (i.e of model category structure) involved, and this operation corresponds to somehow "put the objects in general position". This process makes the objects of derived geometry quite flabby to manipulate with the many advantages it might have but also the drawbacks they drag behind : The few structures that survive this "relaxing of strictness" are often invariants of a specific kind. In the case of derived non-commutative schemes, a major class among them are Morita invariants, these are precisely the functors (or higher functors) whose value is preserved by taking triangulated hulls. \footnote{Recall that one of the main tasks of dg-category theory is to give a functorial enhancement for the cone operation in a triangulated category, so the notion of pretriangulated and triangulated dg-category is central.} In this paper we will show how to apply the derived Morita theory on dg-categories to convert theorems concerning Morita invariants of dg-algebras toward theorems concerning non-commutative schemes, by "converting" we really mean that we will prove exportation theorems stating that if an identity involving a Morita invariant holds on some class of dg-algebras satisfying a property $P$ in the underived setting, then this relation persists in the derived setting under mild conditions regarding the stability by cofibrant replacement of the objects satisfying the property $P$. If one is trying to export an identity involving a single Morita invariant -- e.g the Hochschild cohomology $\Hoc^{\star}(A)$ -- evaluated at a single non-commutative scheme $\xX=\perf(A)$, the exportation procedure might sometimes be trivial. For instance, as $\xX$ is the triangulated hull (in the dg-categorical sense) of the classical algebra $A$, elementarily we know from Morita invariance that $\Hoc^\star(\xX)\cong \Hoc^\star(A)$, and so it is easy to see that the Hochschild cohomology $\Hoc^\star(\xX)$ of a non-commutative scheme bears the same well known structure of Gerstenhaber algebra than its little brother $\Hoc^\star(A)$. However the relationship between Morita invariants of non-commutative schemes and their underlying algebras becomes more sophisticated in situations where Morita invariants interact with the derived symmetric monoidal structure $\otimes^L$ over dg-categories and $\widehat{\otimes}^L$ between triangulated dg-categories. For instance, it is not automatic that a relation such as a Künneth formula :  \[\Hoc^\star(A\otimes B)\cong\Hoc^\star(A)\otimes\Hoc^\star(B)\]

holding for a certain class of dg-algebras, can indeed be converted into a Künneth formula involving a product of derived non-commutative schemes, even with some annoying decorations, because the derived tensor product of non-commutative schemes involve using replacements for some model structure. Our purpose here is to give a generic system of conditions to ensure that such formula can be exported, and show how to systematically convert the corresponding theorem on associative algebras, toward its analogous on non-commutative schemes, both with and without decorations. As we said, the main refinement induced by considering non-commutative schemes compared to associative algebras when proving a theorem about Morita invariants, is that one has to perform cofibrant replacements along the manipulations to prove the desired lemmas, and so one must ensure that the objects and hypotheses involved in proving the formula are stable under cofibrant replacements in a sense that we shall describe. The resulting theory amounts to perform the manipulations over sorts of "algebras in general position" with respect to each other and is part of the general machinery of dg-categories. The main result of the paper is Theorem \ref{GammaKunneth}, whose aim is to extend the domain of validity of the Künneth for any weak\footnote{Weak in the sense that we are here only concerned by the invariance under triangulated hull, of the invariant $\Gamma$ and we don't need any information about its functoriality , which might depends on the choice of invariant :  For instance Hochschild homology has an actual (weak)-2-functoriality, i.e dg-bimodules induces morphisms on the Hochschild homology, but on the contrary Hochschild cohomology only has a restricted (weak)-2-functoriality, only dg-bimodule acting faithfully by tensoring induce equivalences on Hochschild cohomology. See \cite{keller2006differentialgradedcategories,keller:hal-03280332} for a detailed treatment.}  Morita invariant $\Gamma$ (Definition \ref{defweakmo}) satisfying a mild condition of stability by cofibrant replacement.

 \begin{theorem}
     Let $(\cC,\otimes_{\cC})$ be a symmetric monoidal category (not necessarily closed), $\pP,\qQ$ be two subcategories of $\dgcat_R$ and $\Gamma:\dgcat_R\ra \cC$ be a weak Morita invariant. Assume that $\qQ$ is stable by cofibrant replacement for Tabuada's model structure, and that $\Gamma$ satisfies an underived Künneth formula on $\pP\times \qQ$, i.e that for any $\aA\in \pP$ and $\bB\in \qQ$ \[\Gamma(\aA\otimes \bB)\cong \Gamma(\aA)\otimes \Gamma(\bB)\]

     Then, $\Gamma$ satisfies a derived Künneth formula on $\widehat{\pP}\times \widehat{\qQ}$, i.e for any $(\aA,\bB)\in \widehat{\pP}\times \widehat{\qQ}$ holds \[ \Gamma(\aA\widehat{\otimes}^L\bB)\cong \Gamma(\aA)\otimes \Gamma(\bB)\]

\end{theorem}

However, as this result is using a very systematic mechanism, it is superseeded by two more general Theorems ( Theorem \ref{ThDiag} and Theorem \ref{ThExact} ) stating roughly the same thing, but replacing the Künneth formula by any diagram involving Morita invariants, or any exact sequence involving Morita invariants. \\

In typical application situations $\Gamma$ is taken among the classical Morita invariants $\Hoc^\star,\Hoc_\star$, $ \HC^\star,\HC_\star$, $\HP^\star,\HP_\star $, namely Hochschild, cyclic and periodic (co)homologies, and $\pP$, $\qQ$ are sub-dg-categories of the category of associative dg-algebras. Usually, $\qQ$ is  chosen to be a suitable subcategory of algebras that is stable by cofibrant replacement, a choice of $\qQ$ is coined to make the homotopy part in the machinery of the theorem work. On the other side a choice of $\pP$ is typically imposed by the hypothesis required to run the Künneth formula (or any other identity on the invariants) on $\pP$, so $\pP$ encodes the necessary algebraic hypotheses for the Künneth formula. If one sees $\pP$ and $\qQ$ as subcategories of those algebras satisfying some categorical property $P$ and $Q$, then, when $\pP=\qQ$ the result can be understood in the following friendly form :  \\

\begin{theorem} Let $\pP$ be a subcategory of $\dgcat_R$ that is stable by cofibrant replacement for Tabuada's model structure (Theorem \ref{tabuada}). Let $\dD$ be a small category, $\Fun(\dD,\cC)$ the category of $\dD$-shaped diagrams valued in $\cC$. Let $D:\dgcat_R^n\ra \Fun(\dD,\cC)$ be a diagram of multilinear Morita invariants (Definition \ref{multimoritadiag}). Then $D$ is commutative on the full subcategory $(\pP\cap \mathrm{dgalg}_R)^n$ if and only if it is commutative on $(\widehat{\pP}\cap \mathrm{nc}_R)^n$ the category of non-commutative schemes whose underlying dg-algebra satisfies property $\pP$. 
\end{theorem}

The proof of this theorem follows a fairly simple scheme whose complication mainly come from the necessity to take the homotopy theory into account (by deriving tensor products in a model category that is compatible with the Morita invariance theory) :  We must first setup the homotopical model for dg-categories and review its symmetric monoidal structures, which is the objective of Section $1$. Then, in Section $2$ we will prove that the derived tensor product of the symmetric monoidal structure of dg-categories preserves Morita equivalences and finally we should apply this result, in order to deduce our Künneth formula a non-commutative scheme $\xX=\perf(A)$, based on its validity on the underlying dg-algebra $A$. Along the way we discover that our scheme of proof can be extended to many similar situations and we give a general tool box for potential applications. Finally in the last section we give a geometrical application of this theory and we apply our Künneth formula to compute the Hochschild cohomology of a product of geometric quasi-phantoms. This is an opportunity to explore the relationship between the monoidal structures involved in the geometric world (over projective scheme seen through their derived category) and the monoidal structure on the dg-enhancement and to compare different homotopy replacement of objects (especially cofibrant replacements for the Morita model structure and h-projective replacements).

\section{Homotopy and Morita theory for dg-categories}

Before proceeding with the proof of the Künneth formula (Theorem \ref{GammaKunneth}), we need to recover a few ingredients. As explained in introduction the scheme of proof consist in reducing the statement from to the case of dg-categories to the case of dg-algebras (which is well known), using Morita equivalences. Thus in this section we shall recall : 

\begin{enumerate}
    \item The definition of dg-categories, their (raw) monoidal structure.
    \item As Morita invariance only appears at the level of homotopy we shall recall the basic facts about the monoidal structure on the homotopy category of dg-categories.
    \item Then we will define the Hochschild cohomology of a dg-category and recall its Morita invariance.
    \item Finally we shall recall a Morita invariance result concerning the derived tensor product, that will be a key result for the proof.
\end{enumerate} 

\subsection{First definitions}

Here we loosely recollect the elementary definitions of dg-categories, morphism between them, the Yoneda lemma and the definition of the underlying category : 

\begin{definition}\label{dgcatdef}
    
A dg-category $\mathcal{A}$ over a commutative ring $R$ is a category enriched over complexes of $R$-modules $C(R)$, i.e it is the data of : 
\begin{enumerate}
    \item A family of objects $\Ob(\aA)$.
    \item For each pair of objects $x,y\in\Ob(\aA)$ a chain complex of morphisms $\Hom_{\aA}(x,y)$.
    \item A unit $\id_x: R\ra \Hom_{\aA}(x,x)$
    \item A composition defined for any triple 
    $x,y,z$ by a morphism of chain complexes :  $\circ_{x,y,z}:\Hom_{\aA}(x,y)\otimes\Hom_{\aA}(y,z)\ra \Hom_{\aA}(x,z)$.
\end{enumerate}

Satisfying the usual associativity and unity commutative squares.

\end{definition}

\begin{example}The simplest example of a dg-category is $C_{\dg}(R)$ the dg-category of cochain complexes of $R$-modules. For a pair of complexes $x=(x^i,d_x^i),y=(y^i,d_y^i)$,  its complex of morphisms is defined in degree $n$ by the group: \begin{equation}\Hom_{C(R)}^n(x,y)=\prod_{i\in \mathbb{Z}}\Hom_{\text{R-mod}}(x^i,y^{i+n})\end{equation}
with differential given by $d^n:\Hom_{C(R)}^n(x,y)\longmapsto \Hom_{C(R)}^{n+1}(x,y)$, where $d^n(f^i)_{i\in \mathbb{Z}}=d^{i+n}_y\circ f^i+(-1)^{n+1}f^{i+1}\circ d^i_x$. 
\end{example}

\begin{example}
Another elementary dg-category is the category with one single object and complex of morphisms given by a dg-algebra $A$.
\end{example}

  Small dg-categories form a category $\dgcat_R$, whose morphisms are enriched functors over the complex of dg-$R$-modules. 
  The composition of two dg-functors $F:\mathcal{A}\longmapsto \mathcal{B}$ and $G:\mathcal{B}\longmapsto \mathcal{C}$ is a dg-functor. Forgetting the structure of chain complex lead to a classical category of functors and we have a notion of natural transformation. Such natural transformation $\varphi\in \Hom_{\Fun_{dg}(\aA,\bB)}(F,F')$ can be associated naturally a structure of chain complex: it is the subcomplex of $\prod_{x\in \mathcal{A}}\Hom_\mathcal{B}(F(x),F'(x))$ formed by those $x$ that satisfies \begin{equation}F'(f)\circ \varphi_x=(-1)^{\vert\varphi\vert\vert f\vert}\varphi_y\circ F(f).\end{equation}

Then, we can  define the category of dg-categories over $R$ to be the category with objects the class of dg-categories and for morphisms its class of dg-functors. We define the opposite dg-category $\mathcal{A}^\text{op}$ with care about the sign for the composition law  $f\circ_{\mathcal{A}^{op}} g=(-1)^{\vert f\vert\vert g\vert}g\circ_\mathcal{A} f$ (it is often called the Koszul sign). We get a dg-category of dg-presheaves: 

\begin{definition}Let $\mathcal{A}$ be a dg-category. We define a right  dg-$\mathcal{A}$-module as a dg-presheaf $F:\mathcal{A}^{op}\longmapsto C_{\dg}(R)$, and we denote by $C_{\dg}(\mathcal{A})=\Fun_{\dg}(\mathcal{A}^{op},C_{\dg}(R))$ the dg-category of dg-$\mathcal{A}$-module.
\end{definition}

So we can represent any object of any dg-category $\mathcal{A}$ as a particular case of presheaf valued in chain complexes, just as we use to do in classical sheaf theory valued in $\set$. This representation theory arises as a consequence of the categorical enrichment, which upgrades the Yoneda lemma into:

\begin{lemma}[dg-Yoneda]Let $\mathcal{A}$ be a small dg-category. To any object $x\in \mathcal{A}$ we can associate a natural representative dg-$\mathcal{A}$-module $x^\wedge:=\Hom _{\mathcal{A}} (\cdot,x)$ such that: 

\begin{equation}\Big(\mathcal{Y}:x\in\mathcal{A}\longmapsto x^\wedge\in C_{\dg} (R)\Big)\in \mathcal{C}_{\dg}(\mathcal{A})\end{equation}
defines a dg-functor, and for $x\in \mathcal{A}$ and $M\in C_{\dg}(\mathcal{A})$ there exists a natural isomorphism: 
\begin{equation}f\in \Hom_{C_{\dg}(\mathcal{A})}(x^\wedge,M)\longmapsto f(1_x)\in M(x) \end{equation}
\end{lemma}

We will say that a dg-$\aA$-module $M$ is representable, when it is in the image of Yoneda's functor, i.e when there exist an $x\in \aA$ such that $M= x^\wedge $.

\paragraph{The underlying category of a dg-category.}The primary interest of dg-categories is to provide a natural setup palliating the lack of functoriality of cones in triangulated categories and so to obtain a better behaved homotopy theory arising from localizations. In our case, we want to study semiorthogonal decompositions of derived categories of coherent sheaves. Thus we will specifically need the derived and homotopy theoretic constructions naturally associated with dg-categories:\\

To a dg-category $\aA$ one can associate its homotopy category $\Ho(\aA)$, which is the ordinary category, whose objects are the same as $\aA$ and morphisms are obtained from the chain complex of morphisms by applying the $0$-th cohomology functor: \begin{equation}\Hom_{\Ho(\mathcal{A})}(x,y)=\Hh^0\big (\Hom_\mathcal{A}^\star(x,y)\big)\end{equation}

Let us recall the general philosophy in which dg-categories are employed and why they are by nature to be studied through the homotopy theory we present below. Assume we are given a category $\mathcal{S}$ whose objects  $X\in \mathcal{S}$ can be associated a triangulated category $D_X$ in a natural fashion, then the objects of $\sS$ can be considered as algebraic spaces in some sense (e.g the classical case where $X$ are schemes and $X\ra D(\Coh(X))$ defines a derivator). Then,  the hole purpose of dg-category theory is to study the structure of dg-enhancements $\mathcal{D}_X$ of 
$D_X$ i.e dg-categories whose underlying category is $D_X$ and their interrelationships whence mapped through a morphism $X \ra Y$. Especially, when changing triangulated category through a morphism, dg-categories solves the problem of functorially mapping the cone of a morphism. Thus, we are essentially  interested in the structure of dg-categories up to equivalence between their homotopy category. This explains why the homotopy theory recalled below has a pivotal role in dg-category theory. 
\paragraph{Tabuada's model structure on dg-categories}We say that a dg-functor $F:\aA\rightarrow \bB$ is a quasi-equivalence if it induces an equivalence between the homotopy categories $\Ho(F):\Ho(\aA)\rightarrow \Ho(\bB)$ and for any pair of objects $x,y\in \cC$ $\Hom_{\cC}(x,y)$ is a quasi-isomorphism of chain complexes.

\begin{theorem}[Model structure on $\dgcat_R$]\label{tabuada}
The category $\dgcat_R$, admits a structure of cofibrantly generated model category when equiped with the three class: 
\begin{enumerate}
    \item Weak equivalences being quasi-equivalences of dg-categories.
    \item Fibrations being the dg-functors  $F:\mathcal{A}\longmapsto \mathcal{B}$ whose map on morphisms is surjective. And such that for any isomorphism $g:F(x)\longmapsto y$ of $\Ho(\mathcal{B})$ there exist an isomorphism $f$ of $\Ho(\mathcal{A})$ such that $F(f)=g$
\end{enumerate}
\end{theorem}

Concretely the usual axioms $(MC)_k$ are satisfied, one can use the usual homotopy machinery: 

\begin{enumerate}
    \item To associate to $\dgcat_R$ its homotopy category in the sense of Quillen,
    \item  To associate to dg-functors $F:\aA\rightarrow \bB$ a pair of left and right derived functors  forming a Quillen adjunction in the homotopy category.
    \item One can even refine the Quillen homotopy category through simplicial (i.e Dwyer-Kan) localisation to display higher homotopical structures, this is the theory developed by Toën in \cite{toen2006homotopy}.
\end{enumerate}

\begin{definition}
    We define the category $\Hqe$ of homotopy quasi-equivalent classes of $\dgcat_R$ these are obtained as the Quillen homotopy category associated to Tabuada's model structure on $\dgcat_R$. Object of $\Hqe$ are dg-categories and morphism are homotopy classes of dg-functors up to weak equivalence.
\end{definition}

In particular $\Hqe$ is the Gabriel-Zismann localisation of $\dgcat_R$ w.r.t to quasi-equivalences of dg-categories, so its morphisms can be realized as spans $\aA\leftarrow\aA'\rightarrow\bB$ of a quasi-equivalence $W:\aA\leftarrow\aA'$ and a dg-functor $F:\aA'\rightarrow \bB$. We call quasi-functors the morphisms in $\Hqe$ and we denote $\qfun(\aA,\bB)$ for the (external) hom set between $\aA$ and $\bB$.\\

From classical model category theory (see e.g $\cite{Idrissi2022,hovey2007model}$ we know that the class of fibrations $\mathcal{F}$ and weak equivalences $\mathcal{W}$ determines the class of cofibrations $\cC$ as exactly those morphisms having the right lifting property w.r.t to trivial fibrations $\mathcal{W}\cap \mathcal{F}$. In turn this class of cofibration defines a class of cofibrant objects, and a cofibrant replacement functor which is hard to describe in general: it is built out of the argument of small object applied to cellular complexes, so by transfinite limit process. Fortunately for us, this process preserves just enough structure for our purpose and we will only need the following:

\begin{lemma}\label{identityonobject}
    The cofibrant replacement functor $Q$ for the model category structure of Tabuada, Theorem \ref{tabuada}, can be assumed to act as the identity on objects. In particular, the cofibrant replacement of a dg-algebra (seen as a dg-category with one object) is still a dg-algebra. 
\end{lemma}

\proof The proof consist essentially in observing that the construction of the cofibrant replacement functor for the model category of Tabuada \cite{tabuada2004une} doesn't need to add new objects while building a cofibrant the replacement. This has first been observed by Toën in \cite{toen2006homotopy} Proposition 2.3, (2).
\begin{flushright}
    $\square$
\end{flushright}

\paragraph{Pretriangulated dg-categories}

\begin{definition}
   Recall there is a well defined shift functor on $C_{\dg}(R)$ acting on chain complexes $x\in C_{\dg}(R)$ by $x[1]^i:=x^{i+1}$ and on morphisms $f:x\ra y$ by $f[1]:x\ra y[1] $ in the obvious way, increasing the degree by one. This allow, for any dg-category $\aA$ and any object $x\in \aA$ to build $\aA\cup \{x[1]\}$ by adjoining an object $x[1]$ to $\aA$, with the property that for any $y\in \aA$ $\Hom(x[1],y):=\Hom(x,y)[-1]$ and $\Hom(y,x[1])=\Hom(y,x)[1]$. If for any $x\in \aA$ $\aA\cup\{x[1]\}\cong \aA$ we shall say that $\aA$ is closed by shifts. \end{definition}

Given a closed morphism $  f: x \to y  $ of degree 0 in a dg-category $  \cC  $ (i.e., $  df = 0  $), a cone $  \mathrm{Cone}(f)  $ is an object $  z \in \cC  $ such that its Yoneda image $  \widehat{z} = h^z=\Hom_{\cC}(-, z)  $ in the dg-category of dg-modules over $  \cC^{\op}  $ which is isomorphic to the mapping cone of the induced morphism $  \widehat{f}: \widehat{x} \to \widehat{y}  $, where the mapping cone in the dg-category of dg-modules is explicitly the dg-module $  M  $ with $  M= \widehat{y} \oplus \widehat{x}[1]  $ and differential given by  \[d_{\mathrm{Cone}(f)}:= \begin{pmatrix} d_{\widehat{y}} & f \\ 0 &-d_{\widehat{x}} \end{pmatrix} \]
A dg-category $  \cC $ is said to be stable under cones if, for every closed $f:x\ra y$, there exists $z \in \cC$  representing this mapping cone (up to isomorphism in the homotopy category of dg-modules). Using this definition of the cone which lifts the notion of cone in triangulated categories, one can define the pretriangulated hull of a dg-category as the dg-category satisfying the following universal property:

\begin{lemma}\label{ptrhull}Let $\aA$ be a dg-category. Then, there exist a dg-category $\aA^{\text{PreTr}}$ that is closed by shift, cones, direct sums and direct summand, and a fully faithful dg-functor $\aA\rightarrow \aA^{\text{PreTr}}$ satisfying the universal property that any dg functor from $\aA$ to a dg-category $\bB$ closed by shift, cones, direct sums and direct summands factorises uniquely by $\aA\rightarrow \aA^{\text{PreTr}}$. The dg-category $\aA^{\text{PreTr}}$ is called the pretriangulated hull of $\aA$.
\end{lemma}

The proof of this lemma is realized by direct construction of the dg-category satisfying the universal property, it is built as the idempotent completion of the following dg-category:

\begin{definition}
    Let $  \mathcal{C}  $ be a small dg-category. A bounded twisted complex over $  \mathcal{C}  $ is a pair $  (E, q)  $ consisting of:

\begin{enumerate}
    \item A graded object $  E = \bigoplus_{n \in \mathbb{Z}} E^n  $ in $  \mathcal{C}  $ that is bounded, i.e., $  E^n = 0  $ for all but finitely many $  n \in \mathbb{Z}  $, and each $  E^n  $ is a finite direct sum of (shifts of) objects of $  \mathcal{C}  $ in the additive hull $  \mathrm{Add}(\mathcal{C}) $. In other words, $  E  $ belongs to the additive hull of the shifted objects of $  \mathcal{C}  $, with finite total support.
    \item A twisting morphism $q$ of total degree $1$ in the graded endomorphism dg-module of $E$, i.e.,\[q = \sum_{i < j} q_{ij} \quad \text{with} \quad q_{ij} \colon E^i \to E^j\](strictly upper triangular with respect to the grading), such that $  q  $ satisfies the Maurer–Cartan equation\[d(q) + q^2 = 0,\] where $d$ denotes the differential induced componentwise from $  \mathcal{C}  $, and $  q^2  $ is the composition in the graded sense.
\end{enumerate}
\end{definition}

The collection of all bounded twisted complexes over $  \mathcal{C}  $ forms a dg-category, usually denoted $\mathrm{Tw}(\mathcal{C})$. It's morphisms are defined as follow: Let $  (E^\star, q)  $ and $  (F^\star, r)  $ be two bounded twisted complexes (with respective index sets and shifts $  n_i  $, $  m_j  $). A morphism $  f: (E,q) \to (F,r)  $ of degree $  m  $ is a collection of components
$f_{ij} \in \mathcal{C}(E_i, F_j)^{m + m_j - n_i}$
satisfying the compatibility condition
\[d(f) + r \circ f - (-1)^m f \circ q = 0.\]
The differential on the Hom-complex between two twisted complexes is the "twisted" differential:
\[d(f) = d_{\mathrm{naive}}(f) + r \circ f - (-1)^{\deg f} f \circ q,\]
where $d_{\mathrm{naive}}$ is the componentwise differential coming from $\mathcal{C}$. Composition of morphisms is given by matrix multiplication (with the usual signs from the grading).\\

\proof (of Lemma \ref{ptrhull}) First we prove that $\mathrm{Tw}(\cC)$ is stable by shift cones and sums. For shifts: if $  X^\star=(X_i,\alpha_{ij})  $ a twisted complex define $  X^\star[1]  $ by shifting indices $  X_i\mapsto X_{i-1}  $, $  \alpha_{ij}\mapsto(-1)^{j-i}\alpha_{ij}  $ (the sign ensures the Maurer–Cartan equation is preserved). The shift functor is a dg-equivalence. Now for cones, let $  f:X^\star\to  Y^\star  $ be a closed degree-$  0  $ morphism in $  \mathrm{Tw}^b(\aA)  $. The cone $  \mathrm{Cone}(f)  $ is the twisted complex whose underlying graded object is $  X^\star[1]\oplus Y^\star  $ (componentwise direct sum) with twisting data given by the block matrix\[\begin{pmatrix}
d_{X^\star[1]} & f \\
0 & -d_{Y^\star}
\end{pmatrix}\] The Maurer–Cartan equation holds because $  f  $ is closed and both $  X^\star  $ and $  Y^\star  $ satisfy it; the cone is functorial in the usual way. Finally for finite direct sums: The direct sum $  X^\star\oplus Y^\star  $ is taken componentwise on the graded pieces with block-diagonal twisting data. This is clearly a twisted complex and satisfies the universal property of biproducts in the dg-category.\\

To close the proof we need to complete $\mathrm{Tw}(\cC)$ by adding direct summands by passing to the idempotent completion. Objects of $  \mathrm{Idem}(\mathrm{Tw}^b(\aA))  $ are pairs $  (X,e)  $ where $  X\in\mathrm{Tw}^b(\aA)  $ and $  e\in Z^0(\mathrm{End}_{\mathrm{Tw}^b(\aA)}(X))  $ is an idempotent i.e $  e\circ e=e  $. A morphism $  (X,e)\to(X,f)  $ of degree $  q  $ is a morphism $  g:X\to X  $ in $  \mathrm{Tw}^b(\aA)  $ satisfying the intertwining conditions
$g\circ e=g=f\circ g$ with differential inherited from $  \mathrm{Tw}^b(\aA)  $. Composition is the obvious one; the resulting datum is again a dg-category. The inclusion $  \mathrm{Tw}^b(\aA)\hookrightarrow\mathrm{Idem}(\mathrm{Tw}^b(\aA))  $ sending $  X\mapsto(X,\mathrm{id})  $ is fully faithful, hence so is the composite inclusion $  \aA\hookrightarrow\aA^{\mathrm{PreTr}}  $.
$  \aA^{\mathrm{PreTr}}  $ inherits closure under shifts, cones, and finite direct sums from $\mathrm{Tw}^b(\aA)$: each operation applied to a pair $  (X,e)  $ yields a new pair with idempotent $e$ (or a compatible conjugate) preserved by the operation, the direct verification uses the intertwining relations. It is closed under direct summands by construction: given an idempotent endomorphism on $(X,e)$ in $  \aA^{\mathrm{PreTr}}  $, its image under the splitting is represented by a pair $(X,e')$ with the obvious idempotent projection and inclusion satisfying the required relations.
Thus $  \aA^{\mathrm{PreTr}}  $ satisfies all the required closure properties and contains $\aA$ via a fully faithful dg-functor. This closes the proof of lemma \ref{ptrhull}.

\begin{flushright}
    $\square$
\end{flushright}

\begin{definition} A dg-category $\aA$ is said to be pretriangulated if the natural inclusion $\aA\rightarrow \aA^{\text{PreTr}}$ is a quasi-equivalence.
\end{definition}

The universal property of Lemma \ref{ptrhull} produce an inverse of the natural$\aA^{\text{PreTr}}\hookrightarrow(\aA^{\text{PreTr}})^{\text{PreTr}} $, and so $\aA^{\text{PreTr}}$ is naturally pretriangulated. Now notice that $C_{dg}(R)$ is obviously pretriangulated, consequently we deduce that $C_{\dg}(\aA)$ is pretriangulated, so by the universal property of the pretriangulated hull above, we obtain a canonical factorisation of the Yoneda embedding $\mathcal{Y}:\aA\rightarrow C_{\dg}(\aA)$ into $\aA^{\text{PreTr}}\rightarrow C_{\dg}(\aA)$ and this functor is an embedding i.e it is fully faithful, so it is a dg-isomorphism on its image. Then observe (as explained at the beginning of \cite{Orlov_2016} Section 2.3) that the image of this embedding is precisely the sub-dg-category $\mathrm{SF}_{fg}(\aA)$ of $C_{\dg}(\aA)$ of those modules $M$ that are semi-free and finitely generated in the following sense:

\begin{definition}
\begin{enumerate}
    \item We say that a dg-$\aA$-module $N$ is free if it is isomorphic to a direct sum of shifted representable modules, i.e of modules of the form $h^{x_i}[n_i]$ for $x_i\in \aA$, $n_i\in \mathbb{Z}$ and $i$ in some arbitrary set $I$.
    
    \item We say that a module $M$ is semi-free if it there exists a (non-necessarily finite) filtration $0=M_0\subset M_1\subset \cdots $ of $M$ by  whose successive quotients $M_{i+1}/M_i$ are  free dg-$\aA$-module. We denote by $\Sf(\aA)$ the category of semi-free modules. It is a full subcategory of $C_{\dg}(\aA)$ \\

    \item In particular if the filtration can be taken finite (i.e $M_n=M$ for some $n\in \mathbb{N}$) then we say that $M$ is semi-free of finite type (or finitely generated). We denote by $\Sf_{fg}(\aA)$ the category of semi-free finitely generated modules. 
\end{enumerate}
\end{definition}

\begin{lemma}Let $\aA$ be a dg-category, then we have the following equivalent assertions:
\begin{enumerate}
    \item $\aA$ is pretriangulated.
    \item The homotopy category $\Ho(\aA)$ is a triangulated subcategory of $\Ho(C_{\dg}(\aA))$.
\end{enumerate}
    
\end{lemma}

\proof (Sketch) This is a direct consequence of the definition :  $\aA$ is pretriangulated means that $\aA$ is isomorphic to $\aA^{\text{PreTr}}$ the pretriangulated envelope which by construction is the closure of the image of the Yoneda functor by shift sums, and cones, so $\aA$ has a triangulated homotopy category. Reciprocally assume that $\aA$ has a triangulated homotopy category $\Ho(\aA)$, $\aA$  or rather its image in $C_{\dg}(\aA)$ is up to homotopy stable by shift, sums and cones so it is quasi-equivalent to $\aA^{\text{preTr}}$

\paragraph{Dg-enhancements}Now that we have recalled the notion of pretriangulated category we can define precisely the notion of dg-enhancement :

\begin{definition}
    Let $T$ be a triangulated category. A dg-enhancement $\varepsilon:T\rightarrow\Ho(\tT)$ of $T$, is a pretriangulated dg-category $\tT$ together with an exact isomorphism (of triangulated category) $\varepsilon:T\longmapsto \Ho(\tT)$.
\end{definition}

We recall that the existence of dg-enhancements is not granted in full generality as there exists triangulated categories admitting no dg-enhancement at all, however the conditions for their  existence and uniqueness have been vastly studied, notably while enhancing the derived category of quasicoherent sheaves on a noetherian scheme. The main reference for this theory is the article of Canonaco and Stellari \cite{canonaco2018uniqueness}. In section \ref{Cechenhancement}, we shall recall the construction of a practical example of such enhancement for quasi-projective schemes allowing to import sheaf cohomology computations on such scheme to study the dg-category its dg-enhancements.

\paragraph{Enriched model structure on dg-modules}

Let $\aA$ be a dg-category. The category $C_{\dg}(\aA)$  is enriched over the category of complexes of $R$-modules --which one bears a model structure with fibrations  -- and this enrichment is compatible with the model category structure of $C_{\dg}(R)$ in a sense that we recall below, thus in particular the homotopy category of dg-categories is enriched over $D(R)$. Also, we benefit the existence of a model-categorical inner-hom, which comes with specific Quillen adjunction (Lemma \ref{lemm1} and \ref{lemm2}), these will be crucial later on, in order to prove that if $f:\aA\ra \bB$ is a Morita equivalence, then $\id_{\cC}\otimes^L f:\cC\otimes^L\aA\ra \cC\otimes^L\bB$ is again a Morita equivalence (Proposition \ref{tensorMoritainv}).

\begin{definition}[$C_{dg}(R)$-model structure]\label{C(R)model}Let $M$ be a model category. We say that $M$ is a  $C_{dg}(R)$-enriched model category (simply a $C_{dg}(R)$-model category) if we are further given a  tensor product action of $C_{\dg}(R)$ on $M$ through a bifunctor \begin{equation}-\otimes -:C_{\dg}(R)\times M\rightarrow M\end{equation} such that:

\begin{enumerate}
    \item We have natural isomorphism for any $c,c'\in C_{\dg}(R)$, and $m\in M$, \begin{equation}c\otimes (c'\otimes m)\cong (c\otimes_{C_{\dg}(R)}c')\otimes m\end{equation}
    \item For any cofibration $i:c\rightarrow c'$ in $C_{\dg}(R)$ and any cofibration $j:m\to m'$ in $M$, the pushout \begin{equation}c\otimes m'\coprod_{c\otimes m}c'\otimes m\rightarrow c'\otimes m'\end{equation} is a cofibration in $M$.
\end{enumerate}

\end{definition}

That strengthening of the notion of model category is made to ensure that the homotopy category of the model category $M$ is enriched over the homotopy category of $C_{\dg}(R)$ modules.  We refer to the course of Toën \cite{toen:hal-00772841} p.28 for more details. For our purpose, the central example will be $M=C_{\dg}(\aA)$.

\begin{lemma}\label{modelmodule}Let $\aA$ be a dg-category. Then the category of dg-$\aA$-module can be made into a $C_{\dg}(R)$-model category by considering :

\begin{enumerate}
    \item The model category structure induced on $C_{\dg}(\aA)$ by the projective model structure on $C_{\dg}(R)$,  equivalences (resp. fibrations, cofibrations) of dg-$\aA$-modules are morphisms of dg-$\aA$-modules (natural dg-transformation) between $\eta:F\rightarrow G$ such that for any $x\in \aA$, $\eta_x:F(x)\rightarrow G(x)\in \Hom_{C_{\dg}(R)}(F(x),G(x))$ is an equivalence (resp. fibration, resp. cofibration) for the model structure on $C_{\dg}(R)$.

    \item The pointwise action of $C_{\dg}(R)$ on $C_{\dg}(\aA)$ by tensor product: \[\otimes:C_{\dg}(R)\times C_{\dg}(\aA)\ra C_{\dg}(\aA)\] is defined for any $x\in \aA$ by $(c\otimes m)_x:= c\otimes_R m_x$.
\end{enumerate}
\end{lemma}

Now define $\bif(C_{\dg}(R))$ the subcategory of $C_{\dg}(R)$ whose objects are cofibrant complex of $R$-modules for the projective model structure on $C_{\dg}(R)$, i.e bounded below complexes of projective modules. For $M$ a $C(R)$-model category, we define $\bif(M)$ the category of bifibrant objects in $M$. In our case, with $M=C_{\dg}(\aA)$, these are objects $F:\aA\rightarrow C_{\dg}(R)$ such that $0 \rightarrow F$ is a cofibration (every object is fibrant for this model category structure, see \cite{tabuada2004une} for further details).

\begin{definition}The category of derived dg-$\aA$-modules $D(\aA)$ is defined to be the homotopy category of $C_{\dg}(\aA)$ for the model structure given  Lemma \ref{modelmodule}. 
\end{definition}

The category of derived dg-$\aA$-modules $D(\aA)$ is of primary importance, notably while manipulating dg-enhancements over non-smooth algebraic spaces. For instance assume that $\dD_X$ is the Čech dg-enhancement of a complex projective scheme: Then $D(\dD_X)$ is nothing else but $ D(\qcoh(X))$, this is proven in Kuznetsov's \cite{Kuznetsov_2013} section 2.4.

\begin{corollary}
    $D(\aA)$ is enriched over $D(R)$ the derived category of unbounded chain complexes of $R$-modules. 
\end{corollary}

\proof This is a direct application of the fact that the model structure is a $C_{\dg}(R)$-model category structure on $C_{\dg}(\aA)$ and a generic consequence of Hovey's \cite{hovey2007model} enriched model category theory. 

\begin{flushright}
    $\square$
\end{flushright}

\subsection{A world of dg-modules}\label{moduleworld}

In section \ref{Cechenhancement} we will employ Kuznetsov's methods to obtain a very explicit enhancement of the derived category of a projective scheme, that shall be used later on to display a dg-categorical tensor product above $D^b(\Coh(X\times Y))$ in order to compute spectral sequences. In his setting, derived tensor products are defined using h-projective replacements instead of the cofibrant replacement of Tabuada's model category, so we should make explicit the compatibility relations between the two different settings. We take this as an opportunity to clarify many relationships between various subcategories of the dg-category of chain complexes $C_{\dg}(\aA)$ over a dg-category $\aA$ giving a simple diagrammatic picture of the relationships between various analogues of flatness, projectiveness, representability and compactness.\\

The homotopy theory of dg-categories aims to obtain a good description of the category of dg-categories (and also of a fixed dg-category $\aA$), through convenient homotopy models, so it aims to provides functorial (fibrant and cofibrant) replacement, that behave suitably for the computation of topological invariants. However as often in model category theory there are several models available depending on the purpose. \footnote{Many homotopical models for $\aA$ using different and not so obviously compatible replacement functors :  for instance \cite{canonaco2018uniqueness},\cite{canonaco2022internal}, \cite{Kuznetsov_2013}, \cite{Orlov_2016} extensively use the h-projective representations of $\aA$ because it is comfortable to work in this setting while working with enhancements of schemes, whereas the whole Morita theory of $\aA$ developed by Toën \cite{toen:hal-00772841}, \cite{toen2006homotopy} needs cofibrant  replacements of $\aA$ to be worked out. As our result will use both theories, this justifies our need to clarify the relatioship between the different representations}, giving rise to a net of dg-functors relating the various descriptions of $\aA$. As any (quasi-)representation of $\aA$ forcibly factors through the Yoneda embedding, everything takes place in the world of dg-$\aA$-modules : 

\begin{definition}\label{pleindemodules}Let $\aA$ be a dg-category. We define : 
\begin{enumerate}
    \item The dg-category of acyclic dg-$\aA$-modules  $\ac_{\dg}(\aA)$ to be the full sub-dg-category of objects $F\in C_{\dg}(\aA)$ such that for any object $a\in \aA$, $F(a)$ is an acyclic complex in $C_{\dg}(R)$ \emph{(}meaning that $\mathrm{H}^\star(F(a))=0$\emph{)}.
    \item  The dg-category of quasi-representable dg-$\aA$-module $q\rep(\aA)$ to be the sub-dg-category of $C_{\dg}(\aA)$ of modules $M$ such that there exist an object $x\in \aA$ and a natural transformation that is an objectwise quasi-isomorphism $M\ra h^x$, where $h^x$ is the dg-$\aA$-module represented by $x$. \\ 
     \item The dg-category of fibrant and cofibrant dg-$\aA$-modules $\bif(\aA)$ for the $C_{\dg}(R)$-model structure on dg-$\aA$-modules.
    \item The dg-category of perfect dg-$\aA$-modules $\perf(\aA)$ to be the sub-dg-category of $C_{\dg}(\aA)$ of objects $F$ whose projection $\mathrm{H}^0(F)$ in the triangulated category $\Ho(C_{\dg}(\aA))$ is compact (commute with arbitrary sums).
    \item The dg-category of homotopy projective dg-$\aA$-modules $\hproj(\aA)$, or h-projective dg-$\aA$-modules for short, to be the sub-dg-category of $C_{\dg}(\aA)$, whose objects are the $F$ such that for any acyclic dg-$\aA$-module $N$  \[\Hom_{\mathrm{H}^0(C_{\dg}(\aA))}(F,N)=0\]
    \item The dg-category of homotopy flat dg-$\aA$-modules $\hflat(\aA)$ to be the sub-dg-category of $C_{\dg}(\aA)$, whose objects are the $F$ such that for any acyclic dg-$\aA^{op}$-module $N$, $F\otimes_{\aA} N\in C_{\dg}(R)$ is acyclic.
\end{enumerate}
\end{definition}

Now notice the following :  for a fixed dg-category $\aA$ any $h$-projective module is $h$-flat (as asserted by \cite{canonaco2022internal} below their Remark 2.6), any bifibrant dg-$\aA$-module is $h$-projective, any quasi-representable module is bifibrant due to the Yoneda embedding being valued in fibrant and cofibrant complexes and stability of cofibrant and fibrant object by quasi-equivalences. Furthermore from  $\cite{Orlov_2016}$ section 2.3 the inclusion $\aA^{\text{PreTr}}\rightarrow  \Sf_{fg}(\aA)$ is an equivalence. We also know the followings :
\begin{enumerate}
    \item  From Canonaco and Stellari's \cite{canonaco2022internal} (bottom of p.2, also proven in \cite{Keller1994} Section 3.2), that the inclusion $\bif(\aA)\rightarrow \hproj(\aA)$ is a quasi-equivalence
    \item From \cite{Orlov_2016} (top  of p.10) that the embedding $\Sf(\aA)\rightarrow \hproj(\aA)$ is  a quasi-equivalence.  
    \item From \cite{Orlov_2016} that the embedding $\Sf_{fg}(\aA)\rightarrow \perf(\aA)\cap\hproj(\aA)$ is a quasi-equivalence. 
\end{enumerate}

This give rise to the following diagram of inclusions (solid arrows) :  

\begin{equation}\label{bigdiag}\xymatrix@C=0.3cm @R=0.6cm{
& \aA^{\text{PreTr}}\ar[r]^{\cong}& \Sf_{fg}(\aA)\ar[rr]\ar[dr]^{\sim}& &  \Sf(\aA)\ar[ddr]^{\sim}& & C_{\dg}(\aA) 
\\
& &\perf(\aA)\ar@{..>}[d]_{\aA \text{triangulated }\eq \exists}^{\sim} &  \perf\cap \hproj(\aA)\ar[ur]\ar[l] & \\
\aA\ar[rr]^{\sim}\ar[uur]^{\aA\text{ pretriangulated }\eq \sim }& & q\rep(\aA)\ar[r]&\underset{=\widehat{\aA}_{pe}}{\underbrace{\perf\cap \bif}}(\aA)\ar[r]\ar[u]_{\sim}&\bif(\aA)\ar[r]^{\sim}&\hproj(\aA)\ar[r]& \hflat(\aA)\ar[uu]& }\end{equation}

Which further admits an extension into the following square of embeddings proven in the previous section: \[\xymatrix{\mathrm{Tw}^b(\aA)\ar[r]\ar[d] & \mathrm{Tw}(\aA)\ar[d]\\
\mathrm{SF}_{fg}(\aA)\ar[r] & \mathrm{SF}(\aA)}\]

Now something quite wonderful happens if one considers the following (fundamental) property a dg-category might have, which is a strengthening of pretriangulatedness :   

\begin{definition}
       A dg-category $\aA$ is said to be triangulated, if the natural inclusion \[q\rep(\aA)\ra \perf(\aA)\] is a quasi equivalence. 
\end{definition}

In this case, many degeneracies occur, and the inner square of \ref{bigdiag} collapses into dg-isomorphisms :  

\begin{lemma}\label{collapse}If $\aA$ is triangulated, then every perfect dg-module is quasi-representable, thus bifibrant and h-projective, and we have dg-isomorphisms
\begin{equation}\perf(\aA)\cong q\rep(\aA)\cong  \perf(\aA)\cap \bif(\aA)\cong \perf(\aA)\cap \hproj(\aA)\end{equation}
\end{lemma}

\proof The solid arrows of the inner square in the previous diagram is only made of natural inclusions. Thus if $\aA$ is triangulated (i.e $\perf(\aA)$ embeds into $q\rep(\aA)$) the dashed arrow defines an explicit inverse for the solid arrows of the inner square, and in particular all morphisms in the square appears to be dg-isomorphisms.
\begin{flushright}
    $\square$
\end{flushright}

\begin{corollary}
    Every triangulated dg-category $\aA$ is in particular pretriangulated.
\end{corollary}

\proof This is a direct consequence of the $2$-of-$3$ property for the class of quasi-equivalences, applied to the diagram \begin{equation}\xymatrix{& {\aA^{\text{PreTr}}}\ar[d]^\sim\\ 
\aA\ar[r]^\sim \ar[ur]& q\rep(\aA)}\end{equation}

The vertical arrow is the composite  \[\aA^{\text{PreTr}}\overset{\cong}{\ra} \Sf_{fg}(\aA)\overset{\sim}{\ra} \perf\cap\hproj(\aA)\overset{\cong}{\ra} \perf(\aA)\overset{\sim}\ra \perf(\aA)\]

As $\aA$ is triangulated this arrow is a quasi-equivalence by Lemma \ref{collapse}. As the Yoneda map induces a quasi equivalence $\aA\overset{\sim}{\ra}q\rep(\aA)$ and the previous diagram commute, the $2$-of-$3$ property for the class of quasi-equivalence (i.e weak equivalence for Tabuada's model structure)  leads to a quasi-equivalence $\aA\cong \aA^{\text{PreTr}}$, which by definition means that $\aA$ is pretriangulated. 

\begin{flushright}
    $\square$
\end{flushright}

Now assume that $\aA$ is triangulated, if we place ourselves in the category of quasi-isomorphism classes of dg-categories up to homotopy $\Hqe/\text{iso}$, the previous diagram collapses into the chain of inclusions : \\

\begin{equation}\xymatrix@R=0.3cm{\aA\ar@{=}[r]& \perf(\aA)\ar[r]\ar@{=}[d]& \qcoh(\aA)\ar@{=}[d]\ar[r]& \hflat(\aA)\ar[r]&D(\aA)\\ 
& q\rep(\aA)\ar@{=}[d] & \hproj(\aA)\ar@{=}[d]\\
& \Sf_{fg}(\aA) &\Sf(\aA) \ar@{=}[d]\\
&  &  \bif(\aA)}\end{equation}

\begin{remarque}Here we have taken as a definition $\qcoh(\aA):=[\Sf(\aA)]_{iso}$. This is justified by the fact that if $\aA=\mathcal{D}_X$ is a dg-enhancement of a smooth projective scheme then $\Ho(\perf(\aA))=D^b(\Coh(X))$ and $\Ho(\Sf(\aA))=D(\qcoh(X))$, see for instance \cite{Kuznetsov_2013,Orlov_2016}.
\end{remarque}

\begin{remarque}We see in particular that there is a quasi-isomorphism between the category of cofibrant (and fibrant) dg-$\aA$-modules and the category of h-projective dg-$\aA$-modules.  
\end{remarque}

\begin{remarque}Notice also that when $\aA$ is triangulated there is a quasi-equivalence from $\aA$ reaching $\perf(\aA)$, so that in particular $\perf(\aA)\cong \widehat{\aA}_{pe}\cong \aA$
\end{remarque}

\subsection{Model categorical power objects and adjunctions}
Of crucial importance to prove Proposition \ref{tensorMoritainv} -- which is the key technical lemma of this article -- will be the following two adjunction lemmas, stating that any quasi-equivalences between dg-categories (resp. the inclusion of the triangulated hull) induce a Quillen adjunction between the underlying model categories of $M$-valued presheaves, where $M$ is an arbitrary $C_{\dg}(R)$-model category.\\

\begin{definition}\label{coftensqeq}If $M$ is a $C_{\dg}(R)$-model category $M$ we say that it is quasi-flat, when for any cofibrant $x\in \Ob(M)$ and any weak equivalence $\varphi:y\ra z$ the induced morphism: $x\otimes_R \varphi:x\otimes y\ra x\otimes z $ is a weak equivalence. 
\end{definition}

Before proceeding, observe that each dg-category $\aA$ comes with a natural $C_{dg}(R)$-model structure induced from $C_{dg}(\aA)$ by restriction. For two $C_{dg}(R)$-model categories, from general (enriched) model category theory (cf eg \cite{hovey2007model}) there exist a model category of functors from $N$ to $M$, that we denote by $M^{N}$. With these notations Toen states the two following lemmas:

\begin{lemma}[Lemma 1 of \cite{toen:hal-00772841}]\label{lemm1}
    Let $g:\aA\ra \bB$ be a quasi-equivalence between dg-categories, and $M$ be a $C_{\dg}(R)$-model category that is cofibrantly generated and quasi-flat. Then,

\[\begin{tikzcd}
    \Ho(M^{\aA})\arrow[bend left=20, ""{name=A}]{r}{g_!}&\Ho(M^{\bB}) \arrow[bend left=20,""{name=B},below]{l}{g^\star}\arrow[from=A, to=B, symbol=\dashv]{}{}
\end{tikzcd}\]

    forms a Quillen equivalence.
\end{lemma}

\begin{lemma}[Lemma 2 of \cite{toen:hal-00772841}]\label{lemm2}
    Let $\aA$ be a dg-category, $h:\aA\ra \widehat{\aA}_{pe}$ be the natural inclusion in the triangulated hull and $M$ a $C_{\dg}(R)$-model category. Then we have a Quillen equivalence 

    \[\begin{tikzcd}
    \Ho(M^{\aA})\arrow[bend left=20, ""{name=A}]{r}{h_!}&\Ho\big(M^{\widehat{\aA}}\big) \arrow[bend left=20,""{name=B},below]{l}{h^\star}\arrow[from=A, to=B,symbol=\dashv]{}{}
\end{tikzcd}\]
\end{lemma}

\subsection{Monoïdal structure on dg-categories}

Let us denote by $\otimes_R$ the tensor product on the category of $R$-modules, it is a closed symmetric monoidal category. It is a classical fact that $\otimes_R$ induces a closed symmetric monoidal  structure on the category of chains complexes $(C_{dg}(R),\otimes_{C(R)})$, and that it does in turn induce a closed symmetric monoidal structure on the derived category $D(R{-\Mod})$ by left deriving the tensor product functor, i.e by setting $\otimes^L:=-\otimes Q(-)$. This can be done in various manners depending on which model structure one choose on $C(R)$, each leading to a particular choice of cofibrant resolution. The same mechanism generalises to derived dg-categories, and in that case we have a closed symmetric monoidal structure on $\Hqe$:\footnote{In reality as explained in Toën \cite{toen2006homotopy}, this is even a monoidal structure for Dwyer-kan simplicial localisation, but we won't use this property here.}

\begin{theorem}[\cite{toen2006homotopy}]\label{toenmonoidal}Let $Q$ be the cofibrant replacement functor for Tabuada's model structure. The category $\Hqe$, endowed with $-\otimes^L-:=Q(-)\otimes-$, is a closed symmetric monoidal category. Concretely, for $\aA$, $\bB$ and $\cC$ three dg-categories, there exists a dg-category $\RHom(\bB,\cC)$ and functorial isomorphisms :

\begin{equation}\Hom_{\Ho_{Q}(\dgcat_R)}(\aA\otimes^L\bB,\cC)\cong \Hom_{\Ho_{Q}(\dgcat_R)}(\aA, \RHom (\bB,\cC))\end{equation}
\end{theorem}

\begin{lemma}\label{supprimerunQ}There is a quasi-equivalence: \[Q(\aA)\otimes Q(\bB)\cong Q(\aA)\otimes \bB\]
\end{lemma}

\proof Recall from usual model category theory that $q:Q(\bB)\ra \bB$ is a trivial cofibration, and hence in particular a weak equivalence for the model categorical structure (here this means a quasi-equivalence of dg-category). From Lemma \ref{coftensqeq} we hence deduce that 
\[\id_{Q(\aA)}\otimes q:Q(\aA)\otimes Q(\bB)\ra Q(\aA)\otimes \bB\]
is a quasi-equivalence.

\begin{flushright}
$\square$
\end{flushright}

In the previous proof notice in particular, that the appealed lemma rely mostly on the quasi-flatness of the model structure on $\dgcat_R$.

\begin{lemma}\label{modelrhom}Let $\aA$ and $\bB$ be two dg-categories. We have the following identification up to quasi-equivalence of dg-categories: 

\begin{equation}\RHom(\aA,\bB)\cong \rep(\aA,\bB)\cong \bif(C_{\dg}(\aA^{op}\otimes^L \bB))^{\text{rqRep}},\end{equation}

where $\rep(\aA,\bB)$ are dg-functors sending representable dg-$\aA$-modules to representable dg-$\bB$-modules and $\bif(C_{\dg}(\aA\otimes^L \bB^{op}))^{\text{rqRep}}$ is the dg-category of right quasi-representable $\aA$-$\bB$-bimodules, i.e object $M$ of $C_{dg}(\aA\otimes^L \bB^{op})$ such that the dg-functor $M:C_{dg}(\aA)\ra C_{dg}(\bB)$ naturally associated by adjunction factors through representables.
\end{lemma}

We further have a good description of the underlying category $\Ho(\RHom(\aA,\bB))$ by comparing it with $\Ho(\rep(\aA,\bB))$ which is equivalent to the category of right quasi-representable $\aA^{op}\otimes \bB$-modules. Finally, in the proofs of Proposition \ref{tensorMoritainv}, we will need to relate $\bif(C_{\dg}(\aA))$ with $\bif(C_{\dg}(R))$, this is done naturally through the following lemma, which in a literate form sounds like "bifibrant dg-presheaves are dg-presheaves valued in bifibrant complexes (for the implied model category structures)" :

\begin{lemma}[Toën \cite{toen:hal-00772841}]\label{unpointneuf}Let $\aA$ be any dg-category, then 
    \[\bif(C_{\dg}(\aA))\cong \Hom(\aA^{op},\bif(C_{\dg}(R)))\]
\end{lemma}

\subsection{Monoïdal structure on triangulated dg-categories}
The raw derived structure above isn't enough to treat of the Morita theory of dg-categories. The idea is precisely that the monoidal structure on triangulated dg-categories differs from the derived monoidal structure on dg-categories: A Morita equivalence will be a dg-functor inducing quasi-equivalence on the triangulated hull (defined by Proposition \ref{trghull}). Recall that an object $x$ in a category $\aA$ is said to be compact, if the presheaf $h_x$ represented by $x$ commutes with arbitrary coproducts. 

\begin{definition}We say that a dg-category $\tT$ is triangulated, if every compact object in $D(\tT^{op})$ is representable. We denote by $\dgcat^{\Delta}_R$ the category of triangulated dg-categories.
\end{definition}

Of great interest for us is the possibility to attach to any dg-category a triangulated hull :

\begin{proposition}\label{trghull}
    Let $\aA$ be a triangulated dg-category, then the forgetful functor \[\Ho_Q(\dgcat^{\Delta}_R)\rightarrow \Ho_Q(\dgcat_R)\] admits a right adjoint denoted $\widehat{\cdot}$. Equivanlently, every  dg-category $\aA$ naturally admits a triangulated hull. The adjoint functor is given by the Yoneda embedding of $\aA$ in the subcategory $\widehat{\aA}$ of $\bif(C_{\dg}(\aA))$ spanned by compact objects of $\bif(C_{\dg}(\aA))$.
\end{proposition}

 The proof of this proposition is sketched in Toën's lecture $\cite{toen:hal-00772841}$. Though triangulated dg-categories seems to be very abstract objects, they are the cornerstone of Morita invariance theory in derived non-commutative geometry and we need nothing more to define Morita equivalences and invariants in a very broad sense : 




\begin{definition}
    A morphism $\aA\rightarrow \bB$ in $\Hqe$ is called a Morita equivalence if it induces an isomorphism between the triangulated hulls $\widehat{\aA}_{pe}\rightarrow\widehat{\bB}_{pe}$ in $\Ho_Q(\dgcat_R)$.
\end{definition}

\begin{definition}\label{defweakmo}
    Let $\cC$ be a category and $\Gamma:\Ob(\dgcat_R)\ra \Ob(\cC)$ be any function  (e.g $\Gamma=\Hoc$  and $\cC=C_{\dg}(R))$. We shall say that $\Gamma$ is weak Morita invariant (valued in $\cC$), when for any $\aA\in \dgcat_R$ and any morita equivalence $f:\aA\ra \bB$, we have an isomorphism in $\cC$
    \[\Gamma(\aA)\cong \Gamma(\bB).\]
\end{definition}

Also, of interest for the theory -- but somewhat sadly regarding our understanding of algebraic structures encoded by triangulated (dg-)categories -- the monoidal closed structure on dg-categories defined by Theorem \ref{toenmonoidal} isn't compatible with the triangulated hull operation, this gives rise to a different monoidal closed structure on triangulated dg-categories that doesn't share its monoidal tensor product with the one ($\otimes^L$) of  $\Hqe$, but is rather defined by taking its triangulated hull: 

\begin{proposition}\label{symmonoid} The homotopy category $\Ho_Q(\dgcat^{\Delta}_R)$ of triangulated categories admits a closed symmetric monoïdal structure inherited by the monoïdal structure on $\Hqe$. Specifically, it's tensor product functor is defined by $\tT\widehat{\otimes}^L \tT':=\widehat{\tT \otimes^L\tT'}_{pe}$.
\end{proposition}

 The triangulated hull of a dg-category $\aA$ has an interpretation in terms of dg-complexes and of the $C_{\dg}(R)$-model structure presented in the previous section :  it is the subcategory of bifibrant compact objects in $C_{\dg}(\aA)$.

\begin{lemma}\label{whatwewant}
   Let $\aA$ be a dg-category. Then $\widehat{\cC}$ is precisely the full subcategory of compact objects in $\bif(C_{\dg}(\cC))$, denoted by: \[\widehat{\cC}_{pe}\cong \bif(C_{\dg}(\cC))^\omega\]
\end{lemma}

To prove the next proposition we will need to commute opposition of dg-categories with the derived symmetric monoidal structure, hence we need to show that: 

\begin{lemma}
    For $\aA,\cC$ any two dg-categories, then \[\aA^{op}\otimes^L\cC^{op}\cong (\aA\otimes^L \cC)^{op}\]
\end{lemma}

\proof Very concretely we need to show that \[\aA^{op}\otimes Q(\cC^{op})\cong (\aA\otimes Q(\cC))^{op}\]
The fact that $(-)^{op}$ commutes with the underived monoidal structure $\otimes$ is elementary, the non-trivial point is that $Q(\cC^{op})\cong Q(\cC)^{op}$, indeed, $\aA=Q(\cC)$ being cofibrant is equivalent to $\star \ra Q(\cC)$ be a cofibration, which in turn is equivalent to it having the RLP property against all trivial fibrations i.e for any fibration $f:\bB\ra \dD$
\[\begin{tikzcd}
    \star\arrow{r}{}\arrow{d}{}& \bB\arrow{d}{f}\\
    \aA \arrow{r}{}\arrow[dashed]{ru}{} & \dD
\end{tikzcd}\]

Now taking the opposite of this square we equivalently have a filling 

\[\begin{tikzcd}
    \star\arrow{r}{}\arrow{d}{}& \bB^{op}\arrow{d}{f^{op}}\\
    \aA^{op} \arrow{r}{}\arrow[dashed]{ru}{} & \dD^{op}
\end{tikzcd}\]

Recall from Theorem \ref{tabuada} that a trivial fibration is a weak equivalence (i.e a quasi-equivalence of dg-category), which is also a surjective dg-functors such that any isomorphism $\psi:f(x)\ra y$ in $\Ho(\dD)$ has for preimage an isomorphism $\varphi$ in $ \Ho (\bB)$, i.e $f(\varphi)=\psi$. Observe that quasi-equivalence being stable by taking opposite, just as the fibration condition, we deduce that $f:\bB\ra \dD$ is a fibration if and only if $f^{op}:\bB^{op}\ra \cC^{op}$ is. Hence $\aA$ is cofibrant if and only if $\aA^{op}$ is cofibrant consequently the quasi-equivalence $Q(\cC)^{op}\ra \cC^{op}$ (obtain aas opposite of the cofibrant replacement of $\cC$) defines a cofibrant replacement of $\cC^{op}$ now we can conclude by unicity of the cofibrant replacement functor.

\begin{flushright}
    $\square$
\end{flushright}

\begin{proposition}\label{tensorMoritainv}
    Let $f:\aA\rightarrow\bB$ be a Morita equivalence and $\cC$ be a dg-category. Then the induced morphism $\id_{\cC}\otimes f:\cC\otimes^L\aA\rightarrow\cC\otimes^L\bB$
    is a Morita equivalence.
\end{proposition}

\proof We want to show that $\id_{\cC}\widehat{\otimes}^L f:\cC\widehat{\otimes}^L\aA\rightarrow\cC\widehat{\otimes}^L\bB$ becomes an equivalence in $\Ho(\dgcat_R)$, i.e it is a quasi-equivalence. First we translate the problem in terms of an equivalence between $C_{\dg}(R)$-model categories. Using Lemma \ref{unpointneuf} and the monoidal structure on $\Hqe$ we obtain a chain of equivalences 

\begin{equation}\label{chain1}\begin{aligned}
    \bif(C_{\dg}(\aA\otimes^L\cC))& \cong \Hom\big((\aA\otimes^L\cC)^{op},\bif(C_{\dg}(R))\big)\\
    & \cong \Hom\big(\aA^{op}\otimes^L\cC^{op},\bif(C_{\dg}(R))\big)\\
    &\cong \Hom\Big(\aA^{op},\Hom\big(\cC^{op},\bif(C_{\dg}(R))\big)\Big)\\
    & \cong \Hom(\aA^{op},\bif(C_{\dg}(\cC)))\\
    & \cong \bif(C_{\dg}(\cC))^{\aA^{op}}\end{aligned}\end{equation} and similarly for $\bB$. Then using Lemma \ref{lemm2} we get an equivalence \[\Ho(\bif(C_{\dg}(\cC))^{\aA^{op}})\cong \Ho(\bif(C_{\dg}(\cC))^{\widehat{\aA^{op}}})\]

where the last equality comes from the compatibility of the model structure with the enrichment.  As $\widehat{f}^{op}:\aA^{op}\ra \bB^{op}$ is a quasi-equivalence and $\bif(C_{\dg}(\cC))$ is cofibrantly generated and quasi-flat, in the sense that for any cofibrant object $x$ in $M$, and any quasi-isomorphism $e \ra e'$ in $C_{dg}(k)$,
the induced morphism $e \otimes x \ra e' \otimes x$ is an equivalence. We can apply Lemma \ref{lemm1}. Taking $g:=\widehat{f}^{op}$ in argument of the lemma and $M=\bif(C_{\dg}(\cC))$ we obtain a homotopy equivalence \begin{equation}\label{ici}
\Ho(M^{\widehat{\aA}^{op}_{pe}})\cong \Ho(M^{\widehat{\bB}^{op}_{pe}}) \end{equation}
Combining it with (\ref{chain1}) we get that \begin{equation}\label{la}\Ho(\bif(C_{\dg}(\aA\otimes^L\cC)))\cong \Ho(\bif(C_{\dg}(\bB\otimes^L\cC)))\end{equation} 

Now from Lemma \ref{whatwewant}, we know that $\aA\widehat{\otimes}^L\cC$ is precisely the subcategory of $\bif(C_{\dg}(\aA\otimes^L\cC))$ of compact objects (see $\cite{toen:hal-00772841}$ below proposition $6$)\footnote{Beware if you use Toën's lecture on dg-categories as a companion to check our proofs, that $C_{\dg}(\cC)$ in our notations is $\cC^{op}$-$\Mod$ in his notation, anyway these are dg-presheaves over $\cC$.} and the same holds replacing $\aA$ by $\bB$. Finally the equivalence of (\ref{ici}) is just $\Ho(f\widehat{\otimes} \id_{\cC})$ written in terms of morphism of power objects between the involved model categories and (\ref{la}) is just unpacking it.  This closes the proof of proposition \ref{tensorMoritainv}. 

\begin{flushright}
    $\square$
\end{flushright}

To close this section we give a proof that the notion of cofibrant replacement from

\subsection{Noncommutative schemes}

The homotopy theory and derived Morita theory of dg-categories provide a powerful framework in which many classical notions from algebra and geometry find natural higher-categorical refinements. In particular, they suggest a way to extend the notion of “scheme” itself beyond the commutative world.\\

In classical algebraic geometry, a (quasi-compact quasi-separated) scheme $X$ can be recovered, up to isomorphism, from its derived category of perfect complexes $\perf(X)$, or more precisely from a suitable dg-enhancement of it. By the work of Bondal–van den Bergh and subsequent refinements, when $X$ is separated of finite type, $\perf(X)$ admits a dg-enhancement that is quasi-equivalent to $\perf(A) $ for a suitable dg-algebra $A$ cohomologically bounded, and often smooth when $X$ itself is smooth. This observation motivates the following point of view, developed in various forms in derived and noncommutative algebraic geometry (see in particular works of Toën \cite{toen2006homotopy,toen2007moduli,toen:hal-00772841,toen:hal-00772841}, Orlov \cite{orlov2008remarksgeneratorsdimensionstriangulated,Orlov_2016,Orlov_2003,ORLOV2020107096}, and Kontsevitch \cite{Kontsevitch,kontsevich2024notesainfinityalgebrasainfinity},): dg-categories of perfect modules over dg-algebras should be regarded as the noncommutative analogues of schemes. More precisely, an arbitrary (small) pretriangulated dg-category that behaves like the dg-category of perfect complexes on some “space” is called a noncommutative scheme. Such objects naturally inherit notions of smoothness, properness, finite-dimensionality, etc., directly from the rich homotopy theory of dg-categories via their underlying triangulated homotopy categories and bimodule structures. This perspective unifies several phenomena:\\
\begin{enumerate}
\item Ordinary commutative schemes embed into this world via their dg-enhancements of $  \perf(X)  $.\\

\item Noncommutative algebras (or dg-algebras) give rise to “affine noncommutative schemes” $  \perf(A)  $, generalizing $  \mathrm{Spec} A  $ in the derived sense.\\
\end{enumerate}

In particular, in this setting, smooth and proper noncommutative schemes play a role analogous to smooth proper varieties: they are “compact” in a derived sense, admit good duality and trace maps, and form the natural domain for noncommutative motives. In this way, noncommutative schemes provide a bridge between classical algebraic geometry, noncommutative algebra, and homotopical algebra. They allow one to treat derived categories of coherent sheaves, derived categories of modules over algebras, and many other triangulated categories on the same footing, while retaining the full functoriality of dg-enhancements.

\begin{definition}\label{noncomschdef}
    A noncommutative scheme with coefficient ring $R$ is a dg-category $\mathcal{S}$ of the form  $\mathcal{S}=\perf(A)$ for a cohomologically bounded $R$-algebra $A$.
\end{definition}

Equivalently a non-commutative scheme $\sS$ can be defined as the triangulated hull  $ \widehat{BA}_{pe}$  of the dg-category 
 $BA$ with one object and complex of morphism $A$ as explained in $\cite{toen:hal-00772841}$ or according to our diagram (\ref{moduleworld}).

\begin{definition}
    A non-commutative scheme $\sS$ will be said to be finite dimensional if $A$ assumed to be finite dimensional in the previous definition.
\end{definition}

The various notions of smoothness and properness for a non-commutative scheme are modeled to match the underlying structure of pretriangulated dg-category in the following way: 

\begin{definition}Let $\mathcal{A}$ be a small $R$-linear pretriangulated dg-category. We shall say that $\mathcal{A}$ is: 
\begin{enumerate}
  
    \item Smooth when $\aA$ is perfect viewed as an $\aA^e$-module, i.e as the diagonal bimodule $\Delta_{\aA}:=\Hom_{\aA}:\aA^{op}\otimes\aA\ra C_{dg}(R)$ in   $D(\aA^{op}\otimes^L\aA$). 
      \item Proper when  $\Delta_{\aA}(x,y)$ is compact in $D(k)$ for any $x,y\in \aA$.
\end{enumerate}

Recall that triangulated category $\mathcal{T}$ (over a field $k$) is said to be proper if, for every pair of objects $x,y \in \mathcal{T}$, the graded vector space $  \bigoplus_{n \in \mathbb{Z}} \Hom_{\mathcal{T}}(x, y[n])$ is finite-dimensional over $k$. Now notice the following:

\begin{lemma}When $R$ is a field $k$ and $\aA$ is triangulated, then the condition of being proper is exactly equivalent to the condition of it's underlying (triangulated) category $\Ho(\aA)$ to be proper. For $R$ a general ring, the condition of being proper in the dg-sense is stronger.

\end{lemma}
\proof  We first perform the proof in the case of $R=k$ a field. Let $  x,y\in\aA  $. By the very definition of the homotopy category,
\[\RHom_{\mathcal{A}}(x,y[n])\;=\;\mathrm{H}^n\bigl(\RHom_{\aA}(x,y)\bigr)\]
(up to the usual sign convention on grading, which does not affect dimensions). Thus the graded vector space on the left-hand side of the triangulated properness condition is precisely the graded cohomology
\[\bigoplus_{n\in\mathbb{Z}}\RHom_{\mathcal{A}}(x,y[n])\;\simeq\;\bigoplus_{n\in\mathbb{Z}}\mathrm{H}^n\bigl(\Hom_{\aA}(x,y)\bigr)\]
as graded $  k  $-vector spaces (up to reindexing).\\

It remains to show that a complex $  M  $ of $  k  $-vector spaces is perfect if and only if its total cohomology $  \bigoplus_n \mathrm{H}^n(M)  $ is finite-dimensional as a graded vector space.\\

($  \Rightarrow  $) If $  M  $ is perfect, it is quasi-isomorphic to a bounded complex $  P^\star  $ of finite-dimensional $  k  $-vector spaces. Then $  \mathrm{H}^n(M)\simeq \mathrm{H}^n(P^\star)  $, so only finitely many $  \mathrm{H}^n(M)  $ are nonzero (bounded support) and each is finite-dimensional. Hence $  \dim_k\bigoplus_n \mathrm{H}^n(M)<\infty  $.\\

($  \Leftarrow  $) Conversely, suppose $  \bigoplus_n \mathrm{H}^n(M)  $ has finite total dimension, i.e finitely many nonzero finite-dimensional cohomology groups. Over a field $  k  $, every $  k  $-vector space is projective (and injective), so $  M  $ is homotopy equivalent to its homology complex$$H(M)\;:=\;\bigoplus_n \mathrm{H}^n(M)[-n]$$placed in degree $  n  $ with zero differential. This $  H(M)  $ is itself a bounded complex of finite-dimensional $  k  $-vector spaces (by assumption), hence perfect. Therefore $  M  $ is quasi-isomorphic to a perfect complex, i.e., $  M\in\perf(k)  $. Applying this to $  M=\Hom_{\aA}(x,y)  $ for every pair $  x,y  $ completes the equivalence. \\

\begin{flushright}
    $\square$
\end{flushright}

\begin{corollary}In the case where $R$ a ring of finite global dimension, i.e when every $R$-module admit a finite projective resolution of length bounded by $\mathrm{gl.dim}(R)$, the above result hold. For a general ring $R$, the condition of being proper in the dg-sense is stronger than having an underlying proper triangulated category 
    
\end{corollary}

\proof The direct sense of the previous equivalence didn't required $R$ any assumption on $R$ so the second statement of the corollary hold. Now let us prove the converse direction in the case where $R$ is of finite global dimension: \\

($  \Leftarrow  $) Conversely, suppose the total cohomology of $  M  $ is bounded and each $  \mathrm{H}^n(M)  $ is finitely generated over $  R  $. Since $  \operatorname{gl.dim}(R) < \infty  $, every $  R  $-module has a finite projective resolution. In particular, each finitely generated module $  \mathrm{H}^n(M)  $ admits a bounded finite resolution by finitely generated projective modules (of length at most $  \operatorname{gl.dim}(R)  $).One can then construct a bounded complex of finitely generated projective $  R  $-modules that resolves the whole $  M  $: start with a projective resolution of the cohomology and use the finite global dimension to control the lengths so that the resulting complex remains bounded. More formally, the subcategory of $  D(R)  $ consisting of perfect complexes is thick (closed under cones, shifts, and direct summands). The objects with bounded finitely generated cohomology generate this thick subcategory when $  \operatorname{gl.dim}(R) < \infty  $, because one can build a perfect resolution step-by-step using the finite projective dimension.A standard reference fact (which holds precisely when the global dimension is finite) is that a complex $  M \in D(R)  $ with bounded and finitely generated cohomology lies in the thick subcategory generated by $  R  $ (shifted), hence is perfect.Over rings of finite global dimension, the perfect derived category $  \operatorname{Perf}(R)  $ coincides with the full subcategory of $  D(R)  $ consisting of complexes with bounded finitely generated cohomology (up to quasi-isomorphism). This is a direct consequence of the fact that every module has finite projective dimension, so one can truncate and resolve without introducing unbounded terms.\\

Thus, $  M = \operatorname{Hom}_{\mathcal{A}}(x,y)  $ is perfect over $  R  $ if and only if its total cohomology (i.e., the graded Homs in the triangulated category) consists of bounded finitely generated pieces. This holds uniformly for all pairs $  x,y  $, so the two properness conditions are equivalent.\\

\begin{flushright}
    $\square$
\end{flushright}

\end{definition}

From a geometrical perspective, an important example of a non-commutative scheme is the one given by a dg-enhancement of the derived category of coherent sheaves over commutative scheme. The fact that it is indeed a non-commutative scheme justifies the name in the previous definition. We postpone the proof that non-commutative schemes indeed generalises commutative one to the next section because we shall first need the notion of dg-enhancement.

\begin{lemma}
    Let $\tT$ be triangulated category (in the classical sense). Then any admissible sub-category $\aA$ of $\tT$ is triangulated. If furthermore $\tT$ is smooth (resp. proper), then $\aA$ is smooth (resp. proper)
\end{lemma}

This is classical triangulated category theory and we leave the proof to the reader.

\subsection{A few results about dg-algebras}

As we have seen in the previous sections, a non-commutative scheme is by construction controlled by its core dg-algebra, so many proofs of non-commutative geometry are given by explicit manipulations on the core which are then exported by Morita invariance. Hence we shall here present a few canonical computation on dg-algebras, which are of general use to estabish the general theory of dg-categories, these result come from \cite{shklyarov2007serredualitycompacthomologically} and we simply upgrade some of his proofs as a warm-up.\\

Let $A$ be a dg-algebra and $N\in C_{dg}(A)$ a right dg-module
$M\in C_{dg}(A^{op})$ a left dg-module. Then, the tensor
product $N\otimes_k M$ is naturally a right $A^e$-module. Let $P\in C_{dg}(A^e)$. We have an isomorphism 
\[
N\otimes^L_AP\otimes^L_AM\cong (N\otimes_kM)\otimes^L_{A^e}P.
\]

Given $M$ a dg-$k$-module we define its dual $
N^*=\Hom_{C_{dg}(k)}(M, k)
$, if $N$ is furthermore a dg-$A$-module, then we can define $N^\star$ valuewise and $N^*$ inherits a canonical structure of a dg-$A^{op}$-module (i.e a structure of left module, $N\in C_{dg}(A^{op})$. Also, if $M\in C_{dg}(A)$, one can define it's dual as a dg-$A$-module as
\[
N^\vee=\Hom_{C_{dg}(A)}(N, A)
\in C_{dg}(A^{ op})\]. Given a dg-algebra $A$ and $N$ a right dg-$A$-module then $N^{**}\cong N$
on the derived level whenever $N$ has finite dimensional total
cohomology.\\

\begin{lemma}Let $M,N\in C_{dg}(A)$. Both the tensor product $N\otimes_kM^*$
and $A$ are right dg-$A^e$-modules. One has an isomorphism
\[
\Hom_{D(A^e)}(A, N\otimes_kM^*)\cong\Hom_{D(A)}(M, N).
\]
\end{lemma}

\proof We place ourselves in the model category $C_{dg}(A^e)$. For any right $A^e$-module $X$ there is a natural isomorphism of complexes
$$\operatorname{Hom}_{A^e}\bigl(X,\operatorname{Hom}_k(M,N)\bigr)\;\simeq\;\operatorname{Hom}_A(X\otimes_A M,N)$$
given on morphisms by
$$f\;\longmapsto\;\bigl(x\otimes m\mapsto f(x)(m)\bigr).$$
The inverse sends $g$ to the map $x\mapsto(m\mapsto g(x\otimes m))$. Specialising to the diagonal bimodule $X=A$ yields a canonical isomorphism of complexes
$$\operatorname{Hom}_{A^e}\bigl(A,\operatorname{Hom}_k(M,N)\bigr)\;\simeq\;\operatorname{Hom}_A(M,N).$$
Composing with $\theta$ therefore produces a natural morphism of complexes
$$\operatorname{Hom}_{A^e}(A,N\otimes_k M^*)\;\longrightarrow\;\operatorname{Hom}_A(M,N).$$
Now the result follows by considering a cofibrant replacements: Let $P\simeq A$ be a $K$-projective resolution in $C_{dg}(A^e)$ (e.g., the normalised bar resolution). Then
$$\operatorname{Hom}_{D(A^e)}(A,B):=H^0\bigl(\operatorname{Hom}_{A^e}(P,B)\bigr)$$
for any right $A^e$-module $B$. Likewise, if $Q\simeq M$ is a $K$-projective resolution in $C_{dg}(A)$, one has
$$\operatorname{Hom}_{D(A)}(M,N):=H^0\bigl(\operatorname{Hom}_A(Q,N)\bigr).$$
The bar resolution of the bimodule $A$ tensors over $A$ with $M$ to give a $K$-projective resolution of $M$. Consequently the underived adjunction of the previous paragraph, applied to $X=P$, descends to an isomorphism
$$\operatorname{RHom}_{A^e}\bigl(A,\operatorname{Hom}_k(M,N)\bigr)\;\simeq\;\operatorname{RHom}_A(M,N)$$
in $D(k)$. Taking $H^0$ yields
$$\operatorname{Hom}_{D(A^e)}\bigl(A,\operatorname{Hom}_k(M,N)\bigr)\;\simeq\;\operatorname{Hom}_{D(A)}(M,N).$$

\begin{flushright}
    $\square$
\end{flushright}
\begin{proposition}

The functor $(-)^\vee:N\in C_{dg}(A)\mapsto N^\vee\in C_{dg}(A^{op})^{op}$ preserves perfect modules, hence descends to
an equivalence $D_{\perf}(A)\rightarrow
D_{\perf}(A^{op})^{op}$. More
precisely, there is a canonical isomorphism $N\cong(N^\vee)^\vee$.
\end{proposition}

\proof This is due to the fact that the functor $N\mapsto\Hom_{C_{dg}(A)}(N, A)$
preserve twisted modules: First observe that the space
$\Hom_{\mathrm{Tw}(A)}((\oplus_{j=1}^nA[r_j], \alpha), (A,
0))$ consists of $n\times1$-matrices of elements of $A$ (which is isomorphic to $\cong A^{
op}$), with differential and $A$-action exactly of the same form as those in $\mathrm{Tw}(A^{op})$). This leads 
 a functor from $D_{\perf}(A)$ to
$D_{\perf}(A^{{op}})^{{op}}$, because
$\Hom_A(-, A)$ sends direct summands to direct summands.
Furthermore, one has a canonical morphism \[
N\rightarrow\Hom_{A^{{op}}}(\Hom_{A}(N, A), A^{op}).\]
which is an isomorphism when $N$ is a
twisted module. If $N$ is a direct summand of a twisted module $P$,
we have an equivalence \[ P\cong\Hom_{A^{{op}}}(\Hom_{A}(P,
A), A^{op}), \text{ and an embedding }
N\hookrightarrow\Hom_{A^{{op}}}(\Hom_{A}(N,
A), A^{op})
\]
which implies
$N\cong\Hom_{A^{{op}}}(\Hom_{A}(N, A),
A^{op})$. Finally, for an arbitrary perfect module $N$,
$(N^\vee)^\vee\cong\Hom_{A^{{op}}}(\Hom_{A}(pN,
A), A^{op})$, where $pN$ is a
resolution of $N$. This closes the proof of this proposition. 
\begin{flushright}
    $\square$
\end{flushright}

\begin{proposition}
If $N$ is perfect and $M$ is arbitrary then there is a canonical isomorphism
\[
{\RHom}_{A}(N,M)\cong M\otimes^{ L}_A N^\vee.
\]
\end{proposition}
\proof First the assignation 
\[
F: M{\otimes}_A \Hom_{C_{dg}(A)}(N,
A)\rightarrow \Hom_{C_{dg}(A)}(N,M)
\]
given by the evaluation $F(m\otimes f)(n)=mf(n)$, where $m\otimes f\in
M{\otimes}_A \Hom_{C_{dg}(A)}(N, A)$ and $n\in N$ defines a morphism of right dg-$A$-module. Indeed it is immediatly linear by the linearity of $f$ and of the product, furthermore it commutes with differentials: writing $\nabla$ for the differential on the tensor product one has

\[\begin{aligned}
\bigl(d\circ F_{N,M}(m\otimes f)\bigr)(n)
&=d_M(mf(n))-(-1)^{|m|+|f|}mf(d_N(n)),\\
\bigl(F_{N,M}(\nabla(m\otimes f))\bigr)(n)
&=d_M(m)f(n)+(-1)^{|m|}m\bigl((df)(n)\bigr),
\end{aligned}\] and the two expressions coincide by the Leibniz rule in $M$ and the definition of $df$. Passing to cohomology in degree zero (or equivalently, replacing the Hom-complexes by their $H^0$) yields a morphism
$$\overline{F}_{N,M}:M\otimes_A\operatorname{Hom}_A(N,A)\longrightarrow\operatorname{Hom}_A(N,M)$$
in the homotopy category. After deriving with respect to $N$ (i.e., replacing $N$ by a cofibrant resolution, or equivalently replacing the underived tensor product by $\otimes_A^L$ and the underived Hom by $\operatorname{RHom}$) one obtains a natural morphism
$$M\otimes_A^L N^\vee\longrightarrow\operatorname{RHom}_A(N,M)$$
in $D(A)$. It remains to show that this morphism is an isomorphism whenever $N$ is perfect.\\

Suppose first that $N=\bigoplus_{i=1}^r A[n_i]$ is a finite direct sum of shifts of $A$. Then, as graded modules,
$$\operatorname{Hom}(N,A)\simeq\bigoplus_{i=1}^r A[-n_i]$$
(with the obvious left $A$-action) and both $M\otimes_A\operatorname{Hom}(N,A)$ and $\operatorname{Hom}(N,M)$ are finite direct sums of shifts of $M$. The map $F_{N,M}$ is the direct sum of the canonical identifications
$$M\otimes_A A[-n_i]\;\xrightarrow{\;\sim\;}\;M[-n_i]\;\xrightarrow{\;\sim\;}\;\operatorname{Hom}(A[n_i],M).$$
Each of these is an isomorphism of complexes (the first is the usual cancellation $M\otimes_A A\simeq M$, the second is evaluation at $1$). Hence $F_{N,M}$ is an isomorphism of complexes, and so is $\overline{F}_{N,M}$.\\

Now let $N$ be a bounded twisted complex built from finite free modules. Concretely, the underlying graded module of $N$ is a finite free graded module
$$F=\bigoplus_{i=1}^r A[n_i],$$
and the differential of $N$ is of the form $d_F+q$, where $q:F\to F$ is a homogeneous $A$-linear map of degree $1$ satisfying the Maurer–Cartan equation $(d_F+q)^2=0$. (In particular $q$ is strictly upper-triangular with respect to any ordering of the generators compatible with the degrees $n_i$.)
The complex $\operatorname{Hom}(N,A)$ has the same underlying graded module as $\operatorname{Hom}(F,A)$, but its differential is
$$d_{\operatorname{Hom}(F,A)}-\operatorname{ad}_q,\qquad
(\operatorname{ad}_q)(f)=f\circ q$$
(up to the usual Koszul sign). Likewise, $\operatorname{Hom}(N,M)$ differs from $\operatorname{Hom}(F,M)$ by the same inner derivation $\operatorname{ad}_q$. The tensor product $M\otimes_A\operatorname{Hom}(N,A)$ differs from $M\otimes_A\operatorname{Hom}(F,A)$ by $\operatorname{id}_M\otimes\operatorname{ad}_q$.
The underlying graded map of $F_{N,M}$ therefore coincides with the already-established isomorphism
$$F_{F,M}:M\otimes_A\operatorname{Hom}(F,A)\;\xrightarrow{\;\sim\;}\;\operatorname{Hom}(F,M).$$
It remains only to check that this graded isomorphism intertwines the two perturbed differentials. For a pure tensor $m\otimes f$ and a generator $e\in F$ one computes
$$\begin{aligned}
\bigl(d_{\operatorname{Hom}(N,M)}\circ F_{F,M}(m\otimes f)\bigr)(e)
&=d_M\bigl(mf(e)\bigr)-(-1)^{|m|+|f|}mf\bigl((d_F+q)(e)\bigr),\\
\bigl(F_{F,M}\circ(d_{M\otimes\operatorname{Hom}(N,A)})(m\otimes f)\bigr)(e)
&=d_M(m)\,f(e)+(-1)^{|m|}m\bigl((d_{\operatorname{Hom}(F,A)}f)(e)-f(q(e))\bigr).
\end{aligned}$$
The two sides are identical by the Leibniz rule in $M$ and the definition of the differential on $\operatorname{Hom}(F,A)$. Consequently $F_{N,M}$ is an isomorphism of complexes, and so is $\overline{F}_{N,M}$.\\

Finally suppose $N$ is a direct summand of a bounded twisted complex $P$ in the homotopy category of dg-modules (hence also in $D(A)$). Write $i:N\to P$ and $p:P\to N$ for the inclusion and projection, so that $p\circ i=\operatorname{id}_N$. Naturality of evaluation gives a commutative diagram
$$\begin{tikzcd}
M\otimes_A\operatorname{Hom}(P,A)\arrow{d}{\id\otimes i^*}\arrow{r}{\overline{F}_{P,M}}& \operatorname{Hom}(P,M)\arrow{d}{i^*}\\
M\otimes_A\operatorname{Hom}(N,A)\arrow{r}{\overline{F}_{N,M}}& \operatorname{Hom}(N,M)
\end{tikzcd}$$
The top arrow is an isomorphism by the previous paragraph. The image of $M\otimes_A\operatorname{Hom}(N,A)$ under $\overline{F}_{P,M}$ therefore lands in the direct summand $\operatorname{Hom}(N,M)\subset\operatorname{Hom}(P,M)$, and the resulting map is a retract of an isomorphism, hence itself an isomorphism. Deriving yields the required isomorphism
$$\operatorname{RHom}_A(N,M)\;\simeq\;M\otimes_A^L N^\vee$$
in $D(A)$.
The construction is functorial in both variables and independent of the choice of presentation of $N$ as a summand of a twisted complex (any two such presentations are homotopic).

\begin{proposition}\label{prop2}
Let $A,B$ be  dg-algebras. Assume that $N$ is a perfect right dg-$A$-module and $M$ is an arbitrary
right dg-$A\otimes B$-module.
Then there is a natural isomorphism of right dg-$B^{{op}}$-modules
\[
({RHom}_{A}(N,M))^*\cong N\otimes^{ L}_A M^*.
\]
\end{proposition}
\proof Consider the canonical map
\[
G_{N,M}: N\otimes_A M^*\rightarrow (\Hom_{C_{dg}(A)}(N,M))^*
\]
given by $G_{N,M}(n\otimes \nu)(f)=(-1)^{|n|(|\nu|+|f|)}\nu(f(n))$,
where $n\otimes \nu\in N\otimes_A M^*$ and
$f\in\Hom_{C_{dg}(A)}(N,M)$. If $N$ is a direct
summand of a twisted $A$-module then $G_{N,M}$ is
obviously bijective and hence an isomorphism.

\begin{flushright}
    $\square$
\end{flushright}

\subsection{Morita invariant properties of dg-categories}
Most of the relevant properties of a dg-category (smooth, proper, compactly generated) are Morita invariant, in this section we shall prove this. This allows us notably to transfer the properties of a non-commutative scheme to its core algebra and reciprocally. \\

For all this section, fix $\aA,\bB$ two dg-categories and $f:\aA\ra \bB$ a Morita equivalence, i.e a dg-functor such that the induced functor  $\widehat{f}:\widehat{\aA}_{pe}\ra \widehat{\bB}_{pe}$ between triangulated hulls is a quasi-equivalence of dg-categories. 

\begin{lemma}\label{smoothdg}
    Assume that $\aA$ is smooth i.e that $\Delta_{\aA}:=\Hom_{\aA}(-,-):\aA^{op}\otimes\aA\ra C_{dg}(k)$ is compact in $D(\aA^{e})$. Then $\bB$ is smooth.
\end{lemma}


\begin{lemma}\label{properdg}
    Assume that $\aA$ is proper i.e that for any $x,y\in \aA$,  $\Delta_{\aA}(x,y):=\Hom_{\aA}(x,y)$ is perfect (i.e compact) in $D(k)$. Then $\bB$ is proper. 
\end{lemma}

The proof of these lemma is left to the reader.

\subsection{Saturated dg-categories}
Of crucial importance for the purpose of this thesis will be the notion of saturatedness of a category. This notion parallels the notion of being both smooth and proper for a commutative scheme and is a generalisation, it was extensively studied by Toen and Vaquié in \cite{toen2007moduli} and we shall refer. First let us recall that according to the previous section, the category of triangulated dg-categories has a monoidal structure denoted $\widehat{\otimes}^L$. It's unit is the dg-category $\widehat{1}:=\bif(C_{dg}(1))$ the dg-category of cofibrant complexes over the punctual dg-category $1$. Hence the dual of an object is defined as $\cC^\vee:=\RHom(\cC,\widehat{1})$. We say that a triangulated dg-category is dualizable (implicitly, inside $\dgcat^{\Delta}_k$) when the natural morphism \[\RHom(\cC,\widehat{1})\widehat{\otimes}^L\cC\ra \RHom(\cC,\cC)\] is an isomorphism.

\begin{definition}Let $\cC$ be a dg-category. We say that $\cC$ is saturated, when one of the following equivalent assertions is satisfied:

\begin{enumerate}
    \item $\cC$ is triangulated, smooth and proper. 
    \item $\cC$ is dualizable in $\Ho(\dgcat^{\Delta}_k)$.
    \item $\cC$ is Morita equivalent to a smooth and proper dg-algebra.
\end{enumerate}
\end{definition}

The proof that those points are all equivalent is the first part of \cite{toen2007moduli} and we won't reproduce it here.\\

\begin{definition}    
Let $\cC$ be a pretriangulated dg-category over a commutative ring $k$. We say that $\cC$ is compactly generated if the following two conditions hold:
\begin{enumerate}
\item The homotopy category $  \Ho(\cC)  $ admits all small coproducts (direct sums).
\item There exists a small set $  \mathcal{G}  $ of objects of $  \cC  $ (called a set of compact generators) such that:
\begin{itemize}\item Every object $  G \in \mathcal{G}  $ is compact in $  \Ho(\cC)  $. That is, for any family $  \{X_i\}_{i \in I}  $ of objects in $  \Ho(\cC)  $, the natural map \[\bigoplus_{i \in I} \Hom_{\Ho(\cC)}(G, X_i) \longrightarrow \Hom_{\Ho(\cC)}(G, \bigoplus_{i \in I} X_i)\]is a quasi-isomorphism (equivalently, the functor $  \Hom_{\Ho(\cC)}(G, -)  $ from $  \Ho(\cC)  $ to the derived category of $  k  $-modules commutes with arbitrary coproducts).

\item The set $  \mathcal{G}  $ generates $  \Ho(\cC)  $ in the following strong sense: if $  x \in \Ho(\cC)  $ is an object such that\[\Hom_{\Ho(\cC)}(\Sigma^m G, x) \cong 0 \quad \text{(in the derived category of }k\text{-modules)}\]for all $  G \in \mathcal{G}  $ and all integers $  m \in \mathbb{Z}  $, then $  X \cong 0  $.
\end{itemize}
\end{enumerate}
\end{definition}

More synthetically, $\cC$ is compactly generated if its homotopy category $\Ho(\cC)$ is a compactly generated triangulated category (in the classical sense of Neeman \cite{neemansurvey}).

\begin{lemma}
    Let $\cC$ be a saturated dg-category. Then $\cC$ is in particular compactly generated and there exist a Serre automorphism $S_{\cC}:\cC\ra \cC$, i.e a functor $S_{\cC}:\cC\ra \cC$ together with a functorial isomorphism: \[\RHom_{\cC}(x,y)\cong \RHom(y,S_{\cC}(x))^\vee\]
\end{lemma}

\proof By definition a dg-category $  \cC  $ is saturated if and only if it is Morita equivalent in $  \mathrm{Ho}(\dgcat_k)  $ to $\perf(C)$, the dg-category of perfect (i.e. compact) right dg-modules over a dg-algebra $C$ that is both proper (i.e $C$ is perfect as a complex of $k$-modules) and smooth (i.e. $C$ is perfect as a $C^e$-bimodule, where $  C^e = C^{op}\otimes_k C  $ is the enveloping dg-algebra). We may therefore replace $  \cC  $ by $\perf(C)$ without loss of generality as both smoothness properness and triangulatedness are Morita invariant notions.\\

 Since $\cC$ is saturated it is in particular triangulated, so by Toën’s construction (Lemma the homotopy category $\Ho(\cC)$ is equivalent to the full subcategory of compact objects inside the derived category of (right) $C$-modules. In particular $\Ho(\cC)$ coincides with the thick subcategory of $D(C)$ generated by the image of $C$ in $C_{dg}(C)$ (i.e the free module of rank one). Thus $  \Ho(\cC)  $ is compactly generated by a single compact object (the class of $C$). \\
 
 Now let us prove the existence of $  S_{\cC}  $. First $\cC$ is smooth and proper by assumption, so according to Lemma \ref{smoothdg} and Lemma \ref{properdg} we know that $C$ is also smooth and proper, Now we construct the Serre functor explcitly: let $C^* = \RHom_k(C,k)  $ be the $k$-linear dual of $C$ viewed as a complex of $k$-modules. It carries a canonical $C$-$C$-bimodule structure via the action
\[(a'\otimes a'')\cdot a = (-1)^{|a'|(|a|+|a''|)}(a''\cdot a\cdot a')\qquad\text{for }a,a',a''\in C,\]
Define a functor on $\perf(C)$ by derived tensor product with this bimodule:
\[S(N)\coloneqq N\otimes^L_C C^*.\]
It is well-defined between perfect modules because $C$ is smooth and proper. Now we can construct an explicit inverse functor for  $S$ on $  D_{\perf}(C)$, by defining the dualizing bimodule $C^!$ and then observing that $S^{-1}=C^!\otimes^L_C-$, where the dualizing bimodule $C^!\coloneqq\RHom_{C^e}^{op}(C,(C^e)^{op})$ and \[  C^*\otimes^L_C C^!\cong C\cong C^!\otimes^L_C C^*  \]

By Shklyarov's Theorem 4.4. 
Moreover $S$ satisfies the required Serre duality isomorphism functorially: for any two perfect $  C  $-modules $  x,y  $,
\[\RHom_{\cC}(x,y)\cong\RHom_{\cC}(y,S(x))^\vee,\]
where $ (-)^\vee=\RHom_k(-,k) $. This follows from the explicit construction and standard adjunction/duality isomorphisms for perfect complexes (see \cite{shklyarov2007serredualitycompacthomologically} Thm. 4.2 and the appendix, especially Prop. A.3). In particular the mapping complexes $  \RHom_{\cC}(x,y)  $ are perfect (by properness), so the dual is well-defined.\\

The construction $  N\mapsto N\otimes^L_C C^*  $ is realized by a dg-functor on a cofibrant model of $  \perf(C)  $ (tensor product with a fixed bimodule is strictly dg-functorial at the level of cofibrant modules). Thus it lifts to a dg-functor $  S_{\cC}:\cC\to\cC  $ inducing the Serre functor on the homotopy category, together with the required functorial isomorphism of mapping complexes.\\

For a general saturated $\cC$ the statement follows by transporting the above data along a chosen Morita equivalence $\cC\cong\perf(C)$ as it preserves mapping complexes up to quasi-isomorphism and induces an equivalence of homotopy categories. This completes the proof.

\begin{flushright}
    $\square$
\end{flushright}

\subsection{Hochschild cohomology}Let $\mathcal{A}$ be a dg-category. Define the envelopping dg-category of $\mathcal{A}$ by $\mathcal{A}^e:=\mathcal{A}^{op}\otimes \mathcal{A}$.

\begin{definition}We define the Hochschild cohomology of a dg-category $\aA$ as the cohomology of the endomorphism complex of the identity $\id_{\aA}$ of its derived endormorphism dg-category $\RHom(\aA,\aA)$: 
\begin{equation}\Hoc^\star(\aA):=\Hh^\star\big(\RHom(\aA,\aA)(\id_{\aA},\id_{\aA})\big)\end{equation}
\end{definition}

In the previous definition, using the model \[\RHom(\aA,\bB)\cong \bif(C_{\dg}(\aA^{op}\otimes^L \bB))^{\text{rqRep}}\] provided by Lemma \ref{modelrhom}, we can take as definition of the object $\id_{\aA}$ the $\aA\otimes\aA^{op}$-bimodule $\Delta_{\aA}:=\Hom_{\aA}:(-,-)\aA^{op}\otimes\aA\ra C_{dg}(R)$. It is right quasi-representable as by definition $\Hom(-,x)$ is the representing presheaf, moreover, it is cofibrant by virtue of the main diagram of section \ref{moduleworld}. In the sequel, we may abuse notation and write $\Hoc^\star(\aA)$ for the complex \[\RHom(\aA,\aA)(\Delta_{\aA},\Delta_{\aA})\].

\begin{proposition}For any associative algebra $A$ seen as a dg-category with one object and $A$ as complex of morphisms, flat over $R$, its Hochschild cohomology coincide with the more classical definition:
\begin{equation}\Hoc^\star(A)\cong \Ext_{A^e}^\star(A,A)\end{equation}
\end{proposition}

Here $\Ext^\star_{\aA^e}(\aA,\aA)$ denotes the complex of endomorphisms of $\aA$ seen as the diagonal bimodule $\Delta_{\aA}:=\Hom_{\aA}(-,-)$ in $D(\aA\otimes \aA^{op})$.

\proof We have an isomorphism in $\Ho(\dgcat_R)$ \[\RHom(A,A)\cong \bif(C_{\dg}(A^{op}\otimes^LA))^{rqRep}\]

where $\bif(C_{\dg}(A^{op}\otimes^LA))^{rqRep}$ is the full subcategory of the dg-category of cofibrant $A^{op}\otimes A$-bimodules formed by right quasi-representable ones and through this isomorphism $1_A$ is sent to $\Delta_{A}$ the diagonal bimodule. But the derived enrichement of the category of dg-modules in $D(R)$ provides \[\begin{aligned}\Hh^\star(\Hom_{\bif(C_{dg}(A\otimes^LA))^{rqRep}}(\Delta_A,\Delta_{A})) 
&\cong \Hom_{\Ho(\bif(C_{dg}(A\otimes^LA^{op})))}(\Delta_A,\Delta_{A})\\
&=\Ext^\star_{A^{e}}(A,A)\end{aligned}\]

\begin{flushright}
    $\square$
\end{flushright}

\begin{lemma}\label{Moritainvariant}
    The Hochschild cohomology of a dg-category $\aA$, is a weak Morita invariant, i.e for any Morita equivalence $f:\aA\ra \bB$ there is an isomorphism 
    \[\Hoc^\star(\aA)\cong \Hoc^\star(\bB)\]
\end{lemma}

\proof First observe that the Yoneda embedding $\aA\ra \widehat{\aA}:=\bif(C_{dg}(\aA))$ induces an isomorphism of $\Ho(\dgcat_R)$ $\RHom(\widehat{\aA},\widehat{\bB})\cong \RHom(\aA,\widehat{\bB})$ (cf e.g \cite{toen:hal-00772841}, Exercise 34). As we have a natural inclusion $\RHom(\aA,\aA)\ra \RHom(\aA,\widehat{\aA})$, the cospan \[\RHom(\aA,\aA)\rightarrow \RHom(\aA,\widehat{\aA})\leftarrow \RHom(\widehat{\aA},\widehat{\aA})\] defines a morphism of $\Ho(\dgcat_R)$, as it is a composite of an inclusion and a quasi-isomorphism, it is quasi-fully faithful. Hence, we obtain quasi-isomorphisms between the $\Hom$-sets of these dg-categories, and especially \[\Hoc^\star(\aA)\cong \RHom(\aA,\aA)(\id_{\aA},\id_{\aA})\cong \RHom(\widehat{\aA},\widehat{\aA})(\id_{\widehat{\aA}},\id_{\widehat{\aA}})\cong \Hoc^\star(\widehat{\aA}).\] 

Of course, this also holds replacing $\aA$ by $\bB$. Now, as $f:\aA\ra \bB$ is a Morita equivalence, $\widehat{f}:\widehat{\aA}\ra\widehat{\bB}$ is a quasi-equivalence and thus we obtain a span of quasi-equivalences \[\RHom(\widehat{\aA},\widehat{\aA})\overset{\cong}{\ra} \RHom(\widehat{\aA},\widehat{\bB})\overset{\cong}{\leftarrow} \RHom(\widehat{\bB},\widehat{\bB})\]

Finally, as $f$ is a dg-functor it matches the identity of $\aA$ to the identity of $\bB$, so we conclude that 

\[\Hoc^\star(\aA)\cong \Hoc^\star(\bB)\]

\begin{flushright}
    $\square$
\end{flushright}

To define the Hochschild homology of a dg-category, we will need the notion of tensor product of a $\aA$-module by a $\aA^{op}$-module "over the base $\aA$": 

\paragraph{Tensor product over a dg-$\aA$-module}

Let $\aA$ be a dg-category. There is a notion of tensor product of dg-modules over $\aA$ generalising the tensor product of modules over an associative algebra, it can be seen as a tensor contraction between a (contravariant) dg-$\aA$-module $M\in C_{\dg}(\aA)=\Fun_{\dg}(\aA^{op},C_{\dg}(R))$ and a covariant dg-$\aA$-module (i.e a dg-$\aA^{op}$-module) $N\in C_{\dg}(\aA^{op})$. It is defined as follows:

\begin{definition}\label{tensorover} We define the tensor product of $M\in C_{\dg}(\aA^{op})$ and $N\in C_{\dg}(\aA)$ over $\aA$, as the chain complex \begin{equation}M\otimes_{\aA }N:=\coker \Big(\bigoplus_{x,y\in \aA}M(x)\otimes \Hom_{\aA}(x,y)\otimes N(y)\overset{\Sigma}{\rightarrow} \bigoplus_{z\in\aA}M(z)\otimes N(z)\Big)\in C_{\dg}(R)\end{equation}

where $\Sigma(v,f,w)=N(f)(v)\otimes w-(-1)^{pq}v\otimes M(f)(w)$, where $p=\deg v$ $q=\deg w$. This definition straightforwardly generalises into a tensor product over $\bB$ between an $\aA^{op}\otimes \bB$-module $M$ and a $\bB^{op}\otimes \cC$-module $N$, defining a new $\aA^{op}\otimes \cC$-module, given for any $x',y'\in \aA^{op}\otimes\cC$ by:
\end{definition}

\begin{equation}M\otimes_{\aA }N(x',y'):=\coker \Big(\bigoplus_{x,y\in \bB}M(x',x)\otimes \Hom_{\bB}(x,y)\otimes N(y,y')\rightarrow \bigoplus_{z\in\bB}M(x',z)\otimes N(z,y')\Big)\end{equation}

\begin{remarque}One can further generalise the previous product straightforwardly to the case where $N$ is a (right) dg-$\aA$-module valued in a dg-category $\dD$, i.e a dg-functor from $\aA^{op}\ra \dD$. In that case, the expression of the 
\end{remarque}

Using the tensor product over a dg-category and the definition of the diagonal $\aA$-$\aA$ bimodule of a category $\aA$ given by $\Delta_{\aA}:=\Hom_{\aA}(-,-)$ we can now recast the usual definition of the Hochschild cohomology of $\aA$:

\begin{definition}We define the Hochschild homology of a dg-category $\aA$ as the cohomology of the tensor product complex of  the diagonal bimodule $\Delta_{\aA}$ of $\aA$ over $\aA\otimes \aA^{op}$:

\begin{equation}\Hoc_\star(\aA):=\Delta_{\aA}\otimes_{\aA\otimes \aA^{op}}^L\Delta_{\aA}\end{equation}
\end{definition}

This definition of the Hochschild homology coincides in the case of a dg-category with only one object, with the usual Hochschild cohomology of the algebra of endomorphisms of this object. Indeed $\Delta_{\aA}$ in that case is the chain complex of the bar resolution. \\

\section{Künneth theorems for noncommutative schemes}

 Using the Morita invariance tool of the previous section, we can now prove the following theorem:

\begin{theorem}\label{GammaKunneth}
     Let $(\cC,\otimes_{\cC})$ be a symmetric monoidal category (not necessarily closed), $\pP,\qQ$ be two subcategories of $\dgcat_R$ and $\Gamma:\dgcat_R\ra \cC$ be a weak Morita invariant. Assume that $\qQ$ is stable by cofibrant replacement for Tabuada's model structure, and that $\Gamma$ satisfies an underived Künneth formula on $\pP\times \qQ$, i.e that for any $\aA\in \pP$ and $\bB\in \qQ$ \[\Gamma(\aA\otimes \bB)\cong \Gamma(\aA)\otimes \Gamma(\bB)\]

     Then, $\Gamma$ satisfies a Künneth formula on $\widehat{\pP}\times \widehat{\qQ}$, both underived and derived, i.e for any $(\aA,\bB)\in \widehat{\pP}\times \widehat{\qQ}$ holds \[ \Gamma(\aA\widehat{\otimes}^L\bB)\cong \Gamma(\aA)\otimes \Gamma(\bB)\]

\end{theorem}

In this theorem we have been denoting $\widehat{\pP}$ the full subcategory of $\dgcat_R$ containing all dg-categories of the form $\widehat{\aA}$ for $\aA$ in $\pP$.





\proof Let $\aA\in \widehat{\pP}$ and $\bB\in \widehat{\qQ}$. By construction $\aA$ is Morita equivalent to an object $a\in \pP\subset \dgcat_R$, i.e $\widehat{a}\cong \aA$. The same holds true for $\bB$ which is Morita equivalent to some $b\in \qQ$. By assumption we know that $\qQ$ is stable by the cofibrant replacement functor $Q:\dgcat_R\ra \dgcat_R$, and that for any  $a\in \pP$ and $b'\in \qQ$, and by the previous reasonning, we know that the Künneth formula holds in the derived sense, i.e 
\begin{equation}\label{prekun}
    \Gamma(a\widehat{\otimes}^L b')\cong \Gamma(a)\otimes_{\cC}\Gamma(b').
\end{equation} Therefrom we will prove that it also holds on all $\pP\times \qQ$, i.e that $\Gamma(\aA\widehat{\otimes}^L \bB)\cong \Gamma(\aA)\otimes_{\cC}\Gamma(\bB)$. As $\Gamma$ is by assumption a weak Morita invariant, we first obtain $\Gamma(a)\cong \Gamma(\aA)$ and $\Gamma(b)\cong \Gamma(\bB)$, and also that:  

\begin{equation}\Gamma(\aA\widehat{\otimes}^L\bB)\cong \Gamma(\aA\otimes^L \bB). \label{kunnethassumption}\end{equation}.

 From the general properties of model categories we know that the cofibrant replacement functor $Q$ is a weak equivalence for the model structure, here this means that it is a quasi-equivalence of dg-categories, as  quasi-equivalences are in particular Morita equivalences, we further get that $\Gamma(Q(b))\cong \Gamma(b)\cong \Gamma(\bB)$. Now $i_{a}:a\ra \widehat{a}\cong \aA$ is a Morita equivalence, so $b\otimes^L i_{\aA}$ is also a Morita equivalence by Proposition \ref{tensorMoritainv}. Also, the same holds true for the inclusion $i_{b}:b\ra \widehat{b}\cong\bB$, thus \[\Gamma(a\otimes^L b)\cong \Gamma(\aA\otimes^L\bB)\]

Now unfolding the definition of $\otimes^L$ we see that
$\Gamma(a\otimes^L b):=\Gamma(a\otimes Q(b))$, but $b$ belong to $\qQ$ which is stable by the cofibrant replacement functor $Q$, thus $Q(b)\in \qQ$. Setting $b'=Q(b)$ in formula \ref{prekun}, we get that:

\[\Gamma(a\otimes^L b)\cong \Gamma(a)\otimes_{\cC}\Gamma(Q(b))\]

Recasting the previous relations we have proven that: \[\Gamma(\aA\widehat{\otimes}^L\bB)\cong \Gamma(\aA\otimes^L\bB)\cong \Gamma(a\otimes^Lb)\cong \Gamma(a)\otimes_{\cC} \Gamma(Q(b))\cong \Gamma(\aA)\otimes_{\cC}\Gamma(\bB)\]

This closes the proof of Theorem \ref{GammaKunneth}. \begin{flushright}
    $\square$
\end{flushright}

If one wants to get rid of the derived tensor product decoration, one can ask for the following (very restrictive) condition to be satisfied:

\begin{corollary}
    With the notations of the previous lemma: If furthermore $a$ and $b$ are such that $Q(a\otimes b)\cong Q(a)\otimes Q(b)$, then \[\Gamma(a\otimes b)\cong \Gamma(a)\otimes \Gamma(b)\]
\end{corollary}

\proof Recall that $Q$ is a weak equivalence hence in particular a Morita equivalence, so that by weak Morita invariance of $\Gamma$ hold \[\Gamma(a\otimes b)\cong \Gamma (Q(a\otimes b))\]

Now applying the hypothesis leads the left identity below, and the right identity is given by Lemma \ref{supprimerunQ}: 
\[\Gamma(a\otimes b)\cong \Gamma(Q(a)\otimes Q(b))\cong \Gamma(a\otimes^Lb)\]

\begin{flushright}
    $\square$
\end{flushright}

Notice that the conditions of the previous corollary are for instance satisfied when $a,b$ are algebras satisfying the coassociative rigidification condition stated in \cite{peroux2021coalgebrasdwyerkanlocalizationmodel} Definition 3.1. Notice that we have in particular proven the following

\begin{theorem}\label{tensorNC} Let $\xX=\perf(A)$, $\yY=\perf(B)$ be two non-commutative schemes. Let $\zZ=\xX\widehat{\otimes}^L\yY=\perf(\mathfrak{C})$ then \[\mathfrak{C}\cong_{\text{Morita}}A\otimes Q(\mathfrak {B})\]

More generally, for any subcategory $\pP$  of $\dgcat_R$, stable by cofibrant replacement and by the underived tensor product $\otimes_R$, then, the subcategory $\mathrm{nc}^{\pP}_R$ of generalised $\pP$-non-commutative scheme (i.e of dg-categories  $\xX$ morita equivalent to $\perf(p)$ for $p\in \pP$), is a monoidal subcategory of $(\dgcat_R,\widehat{\otimes}^L,\mathbf{1})$ 
    
\end{theorem}

One way to interpret this theorem is as one half of the bridge between underived and derived tensorial noncommutative algebraic geometry: Tensor calculus over non-commutative schemes behaves exactly like tensor calculus over algebras "in general position", where putting an algebra in general position means taking it's cofibrant replacement. As a Morita invariant doesn't feel the operation of putting in general position, their algebraic properties are transfered.  The other half of this bridge to the underived world would be provided by a clear understanding of the relationship between the notion of derived Morita equivalence of algebras and that of classical Morita equivalences of algebras, which is still matter of research.\footnote{Classical Morita equivalence on dg-algebra implies derived Morita equivalence, but the converse isn't clear in general.}



\begin{corollary}\label{kunnethhochcoh}
    Let $\aA$ and $\bB$ be noncommutative schemes over a field $k$, one of them being finite dimensional, then \begin{equation}\Hoc^\star(\aA\otimes^L\bB)\cong \Hoc^\star(\aA)\otimes \Hoc^\star(\bB)\end{equation}
\end{corollary}

\proof We show that the hypothesis of Theorem \ref{GammaKunneth} are fulfilled with $\qQ:=\text{dg-Alg}_k$ the dg-category $k$-algebras seen as a  subcategory of $\dgcat_k$ formed by those dg-categories with one single object, $\pP$ the subcategory of $\qQ$ formed algebras that are finite dimensional, and $\Gamma:=\Hoc^\star$. First, $\Hoc^\star$ is a Morita invariant according to Lemma \ref{Moritainvariant}. Second $\widehat{\qQ}\cong \widehat{\text{dg-Alg}_k}$ is exactly the category of non-commutative schemes, as any non-commutative scheme is of the form $\perf(A)\cong \widehat{A}$ for some cohomologically bounded dg-algebra $A$. Third, $\qQ$ is stable by cofibrant replacement, this is a consequence of Toën's \cite{toen2006homotopy} Lemma 2.3 (2) stating that the cofibrant replacement functor can be taken to be the identity on the objects. Finally, the Künneth formula for dg-algebras of the previous form is exactly the following result due to Le and Zhou \cite{LE20141463}. As I struggled to find a reference for this I signal that it is also stated in \cite{angel2017bvalgebra}: 

\begin{lemma}\cite{LE20141463} Let $k$ be a field and $A$ and $B$ be two $k$-algebras such that
one of them is finite dimensional. Then there is an isomorphism of Gerstenhaber algebras
\[\Hoc^\star(A \otimes^L B)\cong \Hoc^\star(A) \otimes \Hoc^\star(B)\] 
\end{lemma}

This conclude our proof of the Künneth formula for Hochschild cohomology. 
\begin{flushright}
    $\square$
\end{flushright}

Another interesting declination of the previous theorem stands for the Hochschild homology of flat dg-algebras, generalising the Künneth formula of Loday \cite{loday93cyclic} (p.124, Theorem 4.2.5)

\begin{corollary}
Let $\aA$ and $\bB$ be two noncommutative schemes over a ring $R$, where $\aA\cong \perf(A)$. Assume that $A$ and $\Hoc_\star(A)$ are flat over $R$ (e.g $R$ is a field or $\aA$ is a smooth and proper non-commutative scheme, or $\aA$ is a flat ), then there is an isomorphism extending the shuffle product on dg-algebras \begin{equation}\Hoc_\star(\aA\otimes^L\bB)\cong \Hoc_\star(\aA)\otimes \Hoc_\star(\bB)\end{equation}
\end{corollary}

\proof  Again, apply the previous theorem with $\qQ$ the category of dg-algebras and $\pP$ its subcategory of flat algebras with flat Hochschild cohomology.

\section{Extension theorems for Morita invariants}

In this section we further generalise the previous reasoning to provide a generic way of exporting results concerning associative algebras toward non-commutative scheme. The simple idea of these results is the following : If we are given a commutative diagram involving multilinear Morita invariants, or an exact sequence of such invariants, which is known to hold when evaluated in any algebras $(A_i)_{i\in I}$ satisfying some categorical property $(\pP_i)_{i\in I}$ that is preserved by cofibrant replacement, then the same commutative diagram/exact sequence, holds on any family of non-commutative schemes $\aA_i$ whose underlying associative algebras satisfies the properties $\pP_i$.\\

The proof of these is mutatis mutandis the same as the proof the Künneth formula (Theorem \ref{GammaKunneth}) but replacing the identity $\Gamma(\aA\otimes^L\bB)\cong \Gamma(\aA)\otimes \Gamma(\bB)$ by a commutative diagram or an exact sequence or an exact sequence of objects of the form $\Gamma_i(\aA_1\otimes \cdots \otimes \aA_n)$. Though simple, it is an interesting generalisation as it will provide a whole microcosm of corollaries that allow to analyses non-commutative schemes just as we can analyse associatives algebras. Before proceeding, we need a few verbose definitions:

\begin{definition}Fix $(\cC,\otimes)$  an $R$-linear, symmetric monoidal category. We define an $n$-linear Morita invariant to be any map of the form : \[F:(\aA_1,\cdots,\aA_n)\in \Ob(\dgcat_R)^n\ra \widetilde{F}(\aA_1\otimes^L\cdots\otimes^L \aA_n)\]

where $\widetilde{F}:\Ob(\dgcat_R)\ra \Ob(\cC)$ is a Morita invariant. We define the set of multilinear Morita invariants to be the semi-ring of maps generated by direct sums and tensor product of $n$-linear Morita invariant (with varying $n$).
\end{definition}

\begin{example}
    Of course, any Morita invariant can be seen as a multilinear Morita invariant with one argument, or even as a multilinear invariant with several arguments but independant of the later ones. \\
\end{example}

One can display the relationships and properties between different $\cC$-valued weak Morita invariants ($\Gamma_d(\aA_1,\cdots \aA_n))_d$ using diagrams, notice however that those diagrams  need not depend functorially on $\aA_1,\cdots, \aA_n$. 

\begin{definition}\label{multimoritadiag}
    A diagram of multilinear Morita invariants valued in a (non-necessarily dg) category $\cC$ is a map $D:(\aA_1,\cdots,\aA_n)\in \dgcat_R^n\ra D(\aA_1,\cdots,\aA_n)\in \Fun(\dD,\cC)$ such that for any index $d\in \Ob(D)$ of a vertex, $(\aA_1,\cdots,\aA_n)\ra \dD(\aA_1,\cdots,\aA_n)(d)$ defines a multilinear Morita invariant.   
\end{definition}

Now with this definition we can prove the following theorem:

\begin{theorem}\label{ThDiag}
    Assume that a diagram of multilinear Morita invariants \[\dD:\dgcat_R^n\ra D(\aA_1,\cdots,\aA_n)\] is commutative over a poly-disc $\pP_1,\times \cdots \times \pP_n\subset \dgcat_R^n$ with each $\pP_i$ stable by cofibrant replacement. Then, $\dD$ is commutative over $\widehat{\pP}_1\times \cdots \times \widehat{\pP}_n$.
\end{theorem}

\proof The proof is a direct consequence of the following lemma, which is a simple recursion applied to theorem $\ref{tensorNC}$ : 

\begin{lemma}
    Let $\aA_i=\perf(p_i)$, where $p_i\in \pP_i$, and each $\pP_i$ is stable by cofibrant replacement, then \[\aA_1\widehat{\otimes}^L\cdots\widehat{\otimes}^L\aA_n\cong_{\text{Morita}} Q(p_1)\otimes \cdots \otimes Q(p_n)\]
\end{lemma}

In consequence, as $\dD$ is a diagram of multilinear morita invariants :\[\dD(p_1,\cdots,p_n)(d)\cong \dD(\widehat{p}_1,\cdots,\widehat{p}_n)(d)= \dD(\aA_1,\cdots,\aA_n)(d) \]

So if $\dD(p_1,\cdots,p_n)$ is commutative if and only if $\dD(\aA_1,\cdots, \aA_n)$ is commutative. Now for any $(\aA_1,\cdots, \aA_n)\in \widehat{\pP}_1\times \cdots, \widehat{\pP}_n$, $\aA_i=\perf(p_i)=\widehat{p_i}$ for some $p_i\in \pP_i$ which proves the desired result. 
\begin{flushright}
$\square $
\end{flushright}

\begin{corollary}\label{weakening}The hypothesis of the previous theorem can be weakened by only requiring that at least one of the $\pP_i$ appearing in the source of each $k$-linear Morita invariant of $\dD$ is stable by cofibrant replacement.
\end{corollary}

This change wouldn't make the proof any different, but would force us to use a notation making explicitly appearent each multilinear invariant $\Gamma_d(\aA_1,\cdots,\aA_n)$ under the form of a tensor product of $k$-linear factors with varying arities. In the same way we have the following : 

\begin{theorem}\label{ThExact}
    Assume that a diagram of multilinear Morita invariants \[\dD:\dgcat_R^n\ra \Fun(D,\cC)\] is exact in a vertex $d\in \Ob(D)$, over $\pP_1,\times \cdots \times \pP_n\subset \dgcat_R^n$ where all the $\pP_i$ are stable by cofibrant replacement. Then, $\dD$ is exact over $\widehat{\pP}_1\times \cdots \times \widehat{\pP}_n$.
\end{theorem}

The proof of this theorem is left to the reader, it is exactly the same as the previous one replacing "commutative diagram", by "exact sequence". As an application we continue our repertory of Künneth theorems by adding the Künneth exact sequence for cyclic homology : 

\begin{definition}[Cyclic homology]We define the cyclic homology $\HC(\aA)$ of a dg-category $\aA$ through \[\HC(\aA):=\Hh^\star(\Delta_{\aA}\otimes_{\aA\otimes \aA^{op}}^L\Delta_{\aA}\otimes_{R}^L\Lambda)\]

where $\Lambda$ is the dg-$R$-algebra $ R[\varepsilon]/(\varepsilon^2)$ with $\varepsilon$ in degree $-1$.
\end{definition}

This definition, adapted from Keller's \cite{keller:hal-03280332} ensures the compatibility with the usual definition on dg-algebras. To see it easily, observe -- or recall -- that $\Delta_{\aA}\otimes_{\aA\otimes \aA^{op}}\Delta_{\aA}$ defines the Hochschild complex over dg-categories, so the change of base coefficient  $(-)\otimes_R\Lambda$ gives back the usual mixed complex defining cyclic homology.\\

\begin{corollary}(Künneth exact sequence for cyclic homology) Let $\aA$ and $\bB$ be two noncommutative schemes such that $\aA$ and $\Hoc_\star(\aA)$ are h-projective (in the sense that the underlying algebra $A$ and $\Hoc_\star(A)$ are h-projective). Then, there is a long exact sequence :  \[\cdots\ra \HC_n(\aA\otimes \bB)\overset{\Delta}{\longrightarrow} \bigoplus_{p+q=n} \HC_p(\aA)\otimes \HC_q(\bB)\ra \bigoplus_{p+q=n-2}  \HC_p(\aA)\otimes \HC_q(\bB)\ra \HC_{n-1}(\aA\otimes \bB)\ra \cdots \]

Where $\Delta$ extends the canonical associative coproduct on the Hochschild homology of dg-algebras.
\end{corollary}

\proof This is exactly the application of our Corollary \ref{weakening} to Loday's \cite{loday93cyclic} corollary 4.3.12, page 132.

\begin{flushright}
    $\square$
\end{flushright}

Notice that likewise we could obtain Conne's periodicity exact sequence (Loday's \cite{loday93cyclic}  Theorem 2.2.1) but in this case, the use of our Theorem is an overkill, as there is no multilinearity in this exact sequence.

\section{Application :  Product of geometric quasi-phantoms.}

In the section 4 we will apply our Künneth formula to compute the Hochschild cohomology of a geometric phantom category $A\boxtimes B$ built out of two quasi-phantoms using the tensor product in the bounded derived category of coherent sheaves. This is an occasion to show how the theory of non-commutative schemes and dg-enhancements of compactly generated triangulated categories allows to treat geometric objects, such as schemes, exactly as derived categories of associative algebras.  A quasi-phantom category $\aA$ is a full triangulated sub-dg-category of a triangulated dg-category $\tT$, hence admitting in particular a well-defined $K$-theory spectrum\footnote{See e.g \cite{toen:hal-00772841} where the $S$-construction, i.e the construction of a well-defined $K$-theory spectrum yielding a Waldhausen structure on the underlying triangulated category of a triangulated dg-category is evoked, see also \cite{neemansurvey} for a clean theoretical framework outside of the dg-categorical frame.}, whose inclusion $\aA\ra \tT$ admits a left adjoint, and such that the triangulated K-theory $K_\star(\aA)$ is a torsion groups, it is called a geometrical phantom if $\tT=D^b(\Coh(X))$ for some projective scheme $X$. Geometrical quasi-phantoms are of great interest especially because of the following phenomenon :  if the order of the torsion group $K_0(\aA)$ and $K_0(\bB)$ of two geometric quasi-phantom $\aA$ and $\bB$ are relatively prime, then the inner tensor product $\aA\boxtimes \bB\in D^b(\Coh(X\times Y))$, has vanishing $K$-theory but isn't a trivial triangulated category (see \cite{Gorchinskiy_2013}). Our goal in this section is thus to give a formula computing an invariant of a phantom that is typically non-trivial : Hochschild cohomology. To do this, we will need to match the homotopy theory arising from Čech (i.e h-projective) dg-enhancements of a projective scheme, with the Morita homotopy theory of dg-categories, which is more suited to study Morita invariance. Along the way we shall also establish the link between the inner tensor product over derived categories $\boxtimes:D^b(X)\times D^b(Y)\ra D^b(X\times Y)$ and the derived tensor product of their respective dg-enhancements.

\subsection{Čech enhancements and compatibility of derived tensor products.}\label{Cechenhancement}

Here we recall a particular case of dg-enhancement of the derived category of a scheme $X$ called Čech enhancement, the theory we expose is taken from Kuznetsov $\cite{Kuznetsov_2013}$ with only very minor modifications (Kuznetsov works over $X\times X$ whereas we work over $X\times Y$), thus we shall be concise. The idea is simply to perform the construction of the Čech complex between arbitrary locally free sheave $F,G$ and to organise the data into a dg-category, this requires to first work with simplicial open coverings of $X$ and then to switch to the associated normalized chain complex.  \\

Recall that for $\aA$ an abelian category, given a cosimplicial object $A^\star:\Delta \ra \aA $ one can define the normalized chain complex functor $N: \mathrm{Cos}(\aA) \ra C_{dg}(\aA)$:

We set $N(A)^\star$ to be the chain complex whose graduates are given by \begin{equation}(NA)^p:=\bigcap_{i=0}^{p-1}\ker\big(V^{p}\overset{s^i_{p-1}}{\longrightarrow} V^{p-1}\big)\text{ with differential,  } d^p_{NA}:=\sum_{i=0}^{p+1}(-1)^i\partial^{i}_{p+1}\end{equation}.

Now if $\aA$ is furthermore monoidal, the previous functor $N$ is known to define a lax monoidal functor $N$, whose laxator between any two cosimplicial objects $A,B$ in $\aA$ is given by the Alexander-Whitney map \[\Delta_{A,B}:N(A)\otimes N(B)\ra N(A\otimes B)\]

It is defined by \[\Delta_{A,B}(a\otimes b):=(\partial^{p+q}_{p+q}\circ \cdots \circ \partial^{p+1}_{p+1})(a)\otimes (\partial^0_{p+q}\circ \cdots \circ \partial^{0}_{q+1})(b)\]

for $a\in A$ of degree $p$ and $b\in B$ of degree $q$. The Alexander-Whitney map is known to admit a homotopy inverse \[\nabla_{A,B}(a\otimes b):=\sum_{p+q=n}\sum_{\tau\in \mathrm{Shuff}(p,q)}(-1)^{\vert\tau\vert}(s_p^{\tau_{p+1}-1}\circ \cdots\circ s^{\tau_{p+q}-1}_{p+q-1})(a)\otimes (s^{\tau_1-1}_1\circ\cdots\circ s^{\tau_p-1}_{p+q-1})(b),\]

for $a,b$ of degree $n$, $\mathrm{Shuff}(p,q)$ are the $(p,q)$-shuffle permutations, meaning that $i\ra \tau_i$ is such that $\tau_1<\cdots \tau_p$ and $\tau_{p+1}<\cdots < \tau_{p+q}$ and $\vert\tau\vert$ is the parity of $\tau$.\\

Let $X$ be a smooth projective scheme over $k$. Any covering $(U_i)_{i\in I}$ can be made into a simplicial object in the category $\mathrm{Sch}_{/X}$ of schemes over $X$ by setting: \begin{equation}\uU_p=\bigsqcup_{(i_0,\cdots,i_p)\in I^{p+1}}U_{i_0}\cap\cdots\cap U_{i_p}\end{equation}

And defining $\uU:\Delta^{op}\ra \text{Open}(X)$ as $\uU([p]):=\uU_p$, and acting on degeneracy and face maps as follow: For $f_r:[p]\ra [p+1]$ the $r$-th face we set each direct summand of $\uU_{f_r}$ to be the open immersion in $\mathrm{Sch}_{/X}$: \[ \bigcap_{j=0}^p U_{i_j}\ \hookrightarrow\ \bigcap_{j\neq r}U_{i_j}\]
For $s_r:[p+1]\ra [p]$ the $r$-th degeneracy, we set each direct summand of $\uU_{s_r}$ to be the identity: 

$$\bigcap_{j=0}^p U_{i_j}\ \xrightarrow{\simeq}\ \bigcap_{j=0}^{p+1}U_{i'_j},$$
where the new $  (p+1)  $-tuple $  (i'_0,\dots,i'_{p+1})  $ is obtained from $  (i_0,\dots,i_p)  $ by repeating the $  r  $-th index (i.e., inserting a copy of $  i_r  $ in position $  r+1  $). Now each complex $F$ of quasi-coherent sheaves  on $X$ give rise to a functor $\mathrm{Sch}_{/X} \ra C_{dg}(k)$ to the category of complexes of vector spaces by sending an object $  f:Y\to X  $ to the dg-complex of $  k  $-vector spaces
\[\Phi_F(f)\ :=\ \Gamma\bigl(Y,\,f^*F\bigr),\]

Hence a covering give rize to a cosimplicial object in $\mathrm{Vect}_k$
\begin{equation}C^p(X,F):=F(\uU_\star)=\bigoplus_{(i_0,\cdots,i_n)\in I^{p+1}}\Gamma(U_{i_0}\cap \cdots \cap U_{i_n},F)\end{equation}

The totalization by direct sum of the associated normalized complex is called the Cech complex of $F$.
We then denote by $\cC_X$ the cosimplicial category with objects being finite complexes of vector bundles over $X$ and with cosimplicial complexes of morphisms \begin{equation}\Hom_{\cC_X}^\star (F,G):=C^\star(X,F^\vee\otimes G)\end{equation}

For  any such cosimplicial category we can define a dg-category $\dD_X:=N\cC_X$ by taking the normalization of the cosimplicial set of morphisms, as the category obtained from $\cC_X$ by sending objects through the identity and associating to any simplicial complex of morphism it's normalisation:

\[\Hom_{N(\cC_X)}(x,y):=N(\Hom_{\cC_X}(x,y))\]

This   a dg-category, with the composition obtained from the composition of cosimplicial complexes by postcomposing by the Alexander-Whitney map $\Delta:N\big(\Hom_{\cC_X}(x,y)\big)\otimes N\big(\Hom_{\cC_X}(y,z)\big)\rightarrow N\big(\Hom_{\cC_X}(x,y)\otimes \Hom_{\cC_X}(y,z)\big)$, with its natural chain complex structure.

\begin{lemma}[Kuznetsov\cite{Kuznetsov_2013}]\label{cechenh}Let $X$ be a smooth projective scheme and $\uU_\star$ an affine open covering. Then the Čech dg-category $\dD_X:=N(\cC_X)$ defines a dg-enhancement of $X$, i.e an equivalence $\varepsilon_X:\Ho(\dD_X)\rightarrow D^b(\Coh(X))$ of $\Ho(\dgcat_\mathbb{C})$ of triangulated categories. This enhancement extends into an equivalence $\widetilde{\varepsilon}:D(\dD_X)\rightarrow D(\qcoh(X))$ such that the following diagram commutes :  \begin{equation}\xymatrix{\Ho(\dD_X)\ar[r]^{\varepsilon}\ar[d]_h&   D^b(\Coh(X))\ar[d]^{\text{inc}}\\ D(\dD_X)&  D(\qcoh(X))\ar[l]_{\widetilde{\varepsilon}}}\end{equation}
\end{lemma}

\proof (Sketch) The functor $  \varepsilon_X \colon \Ho(\dD_X) \to D^b(\Coh(X))  $ is defined on objects by sending a bounded complex $  F  $ of vector bundles in $  \dD_X  $ to the corresponding object of $  D^b(\Coh(X))  $ viewed via the natural embedding of bounded complexes of locally free sheaves into the derived category. On morphisms it is induced by the identity on the underlying complexes. We show that $\varepsilon_{X}$ is an equivalence, by proving fully faithfulness an essential surjectivity. For fully faithfulness, we show that for any two objects $  F,G \in \dD_X  $,
$$H^0\bigl(\Hom_{\dD_X}(F,G)\bigr) \cong \Ext^0_{D^b(\Coh(X))}(F,G).$$
By definition,
$$\Hom_{\dD_X}(F,G) = N\bigl({C}^\bullet(X, F^\vee \otimes G)\bigr).$$
Since $  X  $ is separated (as it is projective) and $  \{U_\alpha\}  $ is an affine open covering, all finite intersections $  U_{\alpha_0, \cdots, \alpha_p}  $ are affine. For any quasicoherent sheaf $  \mathcal{F}  $ on $  X  $, the higher cohomology groups $  H^q(U_I, \mathcal{F}|_{U_I})  $ vanish for $  q > 0  $ (affine schemes have vanishing higher cohomology for quasicoherent sheaves). By the Leray theorem, the Čech complex $  \mathbf{C}^\bullet(X, \mathcal{F})  $ computes the cohomology:
$$H^p\bigl({C}^*(X, \mathcal{F})\bigr) \cong H^p(X, \mathcal{F}).$$
Here $  F^\vee \otimes G  $ is a bounded complex of locally free (hence quasicoherent) sheaves. The normalized complex $  N({C}^*(X, F^\vee \otimes G))  $ is therefore quasi-isomorphic to the hypercohomology complex $  \mathbb{R}\Gamma(X, F^\vee \otimes G)  $. Passing to cohomology yields
$$H^*\bigl(\Hom_{\dD_X}(F,G)\bigr) \cong \Ext^*_{D^b(\Coh(X))}(F,G).$$
In particular the degree-zero cohomology groups agree, so $  \varepsilon_X  $ induces isomorphisms on all Hom-spaces in the homotopy categories, i.e it is fully faithful.

We now prove essential surjectivity. First, every object of $  D^b(\Coh(X))  $ is isomorphic to a bounded complex of coherent sheaves by definition. Since $  X  $ is smooth of finite Krull dimension, every coherent sheaf $  \mathcal{F}  $ admits a finite locally free resolution
\[0 \to E_n \to \cdots \to E_0 \to \mathcal{F} \to 0\]
(with each $  E_i  $ a vector bundle). Thus any bounded complex of coherent sheaves is quasi-isomorphic in $  D^b(\Coh(X))  $ to a bounded complex of vector bundles. Every such complex appears as an object of $  \dD_X  $. Hence every isomorphism class in $  D^b(\Coh(X))  $ lies in the essential image of $  \varepsilon_X  $.\\

Define a dg-functor $  \bar{\varepsilon} \colon C_{dg}(\qcoh(X)) \to C_{dg}(\dD_X)  $ (complexes of quasicoherent sheaves to dg-modules over $  \dD_X  $) by
$$\bar{\varepsilon}(G)(F) := N\bigl(\mathbf{C}^\bullet(X, F^\vee \otimes G)\bigr), \qquad F \in \dD_X.$$
This is a dg-functor because the Čech construction is functorial. It sends acyclic complexes to acyclic dg-modules (again by the Čech computation of hypercohomology), hence descends to a functor
$$\widetilde{\varepsilon} \colon D(\qcoh(X)) \to D(\dD_X).$$
When $  G  $ lies in $  D^b(\Coh(X))  $ and is represented by a bounded complex of vector bundles, $  \bar{\varepsilon}(G)  $ recovers the representable (Yoneda) dg-module up to quasi-isomorphism, so the diagram

\[\xymatrix{\Ho(\dD_X)\ar[r]^{\varepsilon}\ar[d]_h&   D^b(\Coh(X))\ar[d]^{\text{inc}}\\ D(\dD_X)&  D(\qcoh(X))\ar[l]_{\widetilde{\varepsilon}}}\]
commutes (where $  h  $ is the Yoneda embedding).
To see that $  \widetilde{\varepsilon}  $ is an equivalence, note that its restriction to the full subcategory of perfect complexes recovers the equivalence $  \varepsilon_X  $ already established, and every object of $  D(\qcoh(X))  $ can be resolved (in the Čech sense) by objects coming from perfect complexes




\begin{flushright}
    $\square$
\end{flushright}

Being given a triple $\aA,\bB,\cC$ of dg-categories, $M$ and $\aA$-$\bB$-bimodule and $N$ a $\bB$-$\cC$-bimodule, let us define the homotopy projective derived tensor product, $M\otimes^\mathrm{hProj}_{\bB} N:=P(M)\otimes_{\bB} N$ where $P$ is an homotopy projective resolution of $M$, i.e a dg-functor from \[P:C_{dg}(\aA\otimes \bB^{op})\ra \mathrm{hProj}(\aA\otimes \bB^{op})\]. As $\mathrm{hProj}(\cC)$ is quasi-equivalent as a dg-category to $ \bif(\cC)$ (cf. \cite{canonaco2022internal} and the main diagram of our section \ref{moduleworld}), we deduce that the $h$-projective derived tensor product is equivalent to the usual derived tensor product. Therefrom, we see that the following lemma, which is a direct generalisation of a result of Kuznetsov's \cite{Kuznetsov_2013}, applies to the derived setting of the previous section:

\begin{lemma}\label{kuznetxy}
\begin{enumerate}
    \item The Eilenberg-Zilber map \begin{equation}\nabla: N(\cC_X^{op}\otimes \cC_Y)\rightarrow \dD_X^{op}\otimes \dD_Y\end{equation}
    is a dg-functor and defines a  quasi-equivalence of dg-categories. 

 \item The tensor product over $\dD_X$ of a left $\dD_X$-module $M$ and a right $\dD_X$-module $N$  is computed according to \[M\otimes^{\hproj}_{\dD_X}N\cong \Hh^\star\big(X,\varepsilon(M)^\vee\otimes^{\hproj}\varepsilon(N)\big)\]
 \item Let $Y$ be a projective scheme and $\vV$ and affine open covering. Consider $X\times Y$, $\uU_\star \times \vV_\star$ its induced affine covering (with product indexation), then for any coherent sheaf $F$ on $X$ and $G$ on $Y$ holds :  \[C^\star(X\times Y,F\boxtimes G)\cong C^\star(X,F)\otimes C^\star(Y,G)\] 
 Where $\otimes$ is the cosimplicial tensor product. 

 \item For any quadruples of coherent sheaves $F_1,F_2$ on $X$ and $G_1,G_2$ on $Y$ we have an isomorphism \[\Hom_{\cC_{X\times Y}}(F_1^\vee\boxtimes G_1,F_2^\vee\boxtimes G_2)\cong \Hom_{{\cC}_X^{op}}(F_1,F_2)\otimes \Hom_{{\cC}_{Y}}(G_1,G_2)\]
 which is functorial :  it is compatible with compositions and so gives rise to a fully faithful cosimplicial functor \[\kappa:\cC^{op}_X\otimes \cC_Y\ra \cC_{X\times Y} \]

 \item Consequently the cospan $\dD_X^{op}\otimes \dD_Y\overset{\nabla}{\leftarrow }N(\cC^{op}_X\otimes \cC_Y)\overset{\kappa}{\ra} \dD_{X\times Y}$ is a quasiequivalence, and the map $\mu:D(\dD^{op}_X\otimes\dD_Y)\ra D(X\times Y)$ obtained by composing this cospan with the quasi-inverse of $\varepsilon:D(\dD_{X\times Y})\overset{\sim}{\ra}D(X\times Y)$ is quasi-equivalence and so is its restriction to \[\perf(\dD_X^{op}\otimes \dD_Y) \ra D^b(\Coh(X\times Y))\]

 In particular on representable objects, $\mu$ is given by :  \[\mu(h_x\otimes h^y)=\varepsilon(x)^\vee\boxtimes\varepsilon(y)\]
 
\end{enumerate}  
\end{lemma}

The proof of the successive points of the lemma is exactly the content of the section 2.4 of Kuznetsov's article \cite{Kuznetsov_2013}, but replacing line by line $X\times X$ with $X\times Y$, $\cC_X^{op}\otimes \cC_X$ with $\cC_X^{op}\times \cC_Y$, $D(X\times X)$ with $D(X\times Y)$ etc. The reasoning remain unchanged because Kuznetsov never need the hypothesis that $Y=X$. Now consider two admissible subcategories $A\subset D^b(\Coh(X))$ and $B\subset D^b(\Coh(Y))$ and restrict the two dg-enhancements $\dD_X$ and $\dD_Y$ into dg-enhancements $\aA$ and $\bB$ of $A$ and $B$ respectively. \\

Notice that in his paper, Kuznetsov works with Fourier-Mukai kernels (by this we mean bimodules lying in dg-categories of the form $C_{dg}(\cC^{op}\otimes\dD)$) without deriving the tensor product $\cC^{op}\otimes\dD$, which can seem unnatural as we would expect the tensor product to be derived. This is possible due to thesmoothness of $X$, more classically one has:

\begin{theorem}
    Let $X,Y$ be quasi-compact, separated schemes over $k$ such that one of them is flat over $\spec(k)$. Let $\dD_X$ (resp. $\dD_Y$) be any dg-enhancement of $X$ (resp. $Y$). Then \[\bif(\dD_{X\times_k Y})\cong \bif (\dD^{op}_X\otimes^L \dD_Y)\]
\end{theorem}

\proof We start by converting this quasi-equivalence in the notations of Toen \cite{toen2006homotopy} and prove the result using lemmas of him $ \bif(\dD_{X\times_k Y})$ is the dg-category of quasi-coherent complexes over $X\times_kY$ (hear h-projective or quasi-equivalently semi-free complexes) over $\dD_{X\times_k Y}$ a dg-enhancement of $X\times_k Y$, hence by definition of a dg-enhancement it is quasi-equivalent to \[L_{qcoh}(X\times_k Y):=\bif(C_{dg}(\qcoh(X\times_kY))).\]
Similarly $L_{qcoh}(X):=\bif(C_{dg}(\qcoh(X)))$ and same for $Y$.  Now since our pair of scheme satisfy exactly the appropriate hypothesis we can apply \cite{toen2006homotopy} Theorem 8.9 to get \[L_{qcoh}(X\times_k Y)\cong \RHom_{cont}(L_{\qcoh}(X),L_{\qcoh}(Y)).\]
Finally, by \cite{toen2006homotopy} Corollary 7.6 
\[\begin{aligned}\RHom_{cont}(L_{\qcoh}(X),L_{\qcoh}(Y))&\cong \RHom_{cont}(\bif(\dD_X),\bif(\dD_Y))\\
&\cong \bif(\dD_{X}^{op}\otimes^L_{k} \dD_Y)\end{aligned}\]

\begin{flushright}
    $\square$
\end{flushright}

\begin{corollary}
    Let $X,Y$ be two smooth and proper projective schemes over a field $k$. Let $\dD_X$ be the Cech dg-enhancement. Then \[\perf(\dD_{X}^{op}\otimes\dD_Y)\cong \perf(\dD_X\otimes^L\dD_Y)\]
    And:
    \[D(\dD_{X}^{op}\otimes\dD_Y)\cong D(\dD_X\otimes^L\dD_Y)\]
\end{corollary}

\proof Smooth and proper projective schemes over a field satisfies the hypothesis of theorem A and B \cite{canonaco2018uniqueness} (see the discussion below the statements p.3) hence $D^b(\Coh(X\times_k Y))$ admits a unique dg-enhancement, $\perf(\dD_{X}\otimes\dD_Y)$ is a dg-enhancement of $D^b(\Coh(X\times_k Y))$ by Lemma \ref{kuznetxy} 5)  whereas $\perf(\dD_X\otimes^L\dD_Y)$ is a dg-enhancement by the previous theorem. Hence they coincide.

\subsection{Non-commutative schemes generalizes commutative ones}

With the notion of dg-enhancement presented above we are now ready to prove that a commutative scheme is indeed a particular case of non-commutative one. In order to avoid complications we shall restrict ourself to smooth and proper projective schemes over a field, but the result is expected to hold more generally for general schemes and refined notions of dg-enhancements:

\begin{lemma}Let $X$ be a smooth and proper scheme over a field, then $D^b(\Coh{X})$, then $D^b(\Coh(X))$ defines a smooth and proper noncommutative scheme.
\end{lemma}

\proof Let $  \mathcal{C}  $ be any dg-enhancement of $  D^b(\Coh(X))  $ (e.g., the unique one up to quasi-equivalence by the Lunts–Orlov theorem, or the explicit Čech dg-enhancement $\dD_X  $ of Kuznetsov, which satisfies $  \mathrm{Ho}(\dD_X)\cong D^b(\Coh(X))  $ and extends to an equivalence $  D(\dD_X)\cong D(\mathrm{QCoh}(X))  $). We first recall that $  \mathcal{C}  $ (or equivalently $  \perf(X)  $) is itself a smooth and proper dg-category when $X$ is smooth and proper:\\

\begin{enumerate}
\item Properness of $  \mathcal{C}  $: For any objects $  E,F\in\mathcal{C}  $, the Hom-complex $  \mathcal{C}(E,F)  $ is quasi-isomorphic to $  \mathrm{RHom}_X(E,F)  $. Since $  E,F  $ are perfect and $  X  $ is proper, $  \mathrm{R}\Gamma(X,\mathrm{RHom}_X(E,F))  $ is a perfect complex of $  k  $-modules. Hence each $  \mathcal{C}(E,F)  $ is perfect over $  k  $, so $  \mathcal{C}  $ is proper as a dg-category.
\item Smoothness of $  \mathcal{C}  $: The diagonal bimodule $  \Delta_{\mathcal{C}}  $ corresponds (under the equivalence $  \mathcal{C}\cong\perf(X)  $) to the Fourier–Mukai kernel $  \mathcal{O}_\Delta\in D^b(X\times_k X)  $. Since $  X  $ is smooth, the diagonal embedding $  \Delta:X\hookrightarrow X\times_k X  $ is a regular immersion, so $  \mathcal{O}_\Delta  $ admits a finite locally free resolution on $  X\times_k X  $ by its Koszul complex. Thus $ \mathcal{O}_\Delta $ is perfect, and therefore $  \Delta_{\mathcal{C}}  $ is perfect in $  \mathcal{D}(\mathcal{C}^{\mathrm{op}}\otimes_k\mathcal{C})$, so $\mathcal{C}$ is smooth as a dg-category.
\end{enumerate}

It remains to exhibit a dg-algebra $A$ such that $  \mathcal{C}\cong\perf(A)  $. Since $  X  $ is smooth and proper over a field, by Bondal–Van den Bergh \cite{bondal2002generatorsrepresentabilityfunctorscommutative} there exists a perfect complex $  E\in D^b(\Coh(X))  $ that is a classical generator: the smallest thick triangulated subcategory of $  D^b(\Coh(X))  $ containing $E$ is the whole category. Choose a representative of $E$ in $  \mathcal{C}  $. Let
\[A:=\Hom_{\mathcal{C}}(E,E)\]
be its dg-endomorphism algebra. The dg-functor
$$\perf(A)\longrightarrow\mathcal{C},\qquad M\longmapsto M\otimes_A E$$
(defined by sending the free module $  A  $ to $  E  $ and extending by cones, shifts, and direct summands) is then a quasi-equivalence: it is fully faithful on homotopy categories by the generator property and essentially surjective as the thick closure of $E$ is all of $  \mathrm{Ho}(\mathcal{C})$.\\

Finally, since $  \mathcal{C}\simeq\mathrm{Perf}(A)  $ and smoothness/properness of dg-categories are Morita invariants by Lemma \ref{} and Lemma \ref{}, the dg-algebra $A$ is smooth and proper. This shows that the dg-enhancement of $D^b(\Coh(X))$ is quasi-equivalent to $\mathrm{Perf}(A)$ for a smooth and proper dg-algebra $A$, hence defines a smooth and proper noncommutative scheme.

\begin{flushright}
$\square$
\end{flushright}


\begin{flushright}
$\square$
\end{flushright}

\begin{corollary}
    Let $\aA$ and $\bB$ be two pretriangulated full sub-dg-categories of $\dD_X$ and $\dD_Y$ respectively. Then \[\varepsilon_X(\aA)\boxtimes \varepsilon_Y(\bB)\cong \Ho(\perf(\aA^{op}\otimes \bB))\]
\end{corollary}

\proof Under the identification of Lemma \ref{cechenh} point 5, the subcategory $\varepsilon_X(\aA)\boxtimes \varepsilon_Y(\bB)$ of $D^b(\Coh(X\times Y))$ corresponds to the sub-dg-category $\widehat{\aA^{op}\otimes \bB}$ of $\widehat{\dD_X^{op}\otimes \dD_Y}$ indeed :  $\aA^{op}\otimes \bB$ embeds in $\dD_X^{op}\otimes \dD_Y$. Taking the triangulated hull we obtain an inclusion \[\perf(\aA^{op}\otimes \bB):=\widehat{\aA^{op}\otimes \bB}\hookrightarrow \widehat{\dD_X\otimes \dD_Y}:=\perf(\dD_X^{op}\otimes \dD_Y)\]

where the right-hand side is the Čech dg-enhancement of $  D^b(\Coh(X\times Y))  $. To obtain a quasi-equivalence (i.e., the desired isomorphism of dg-categories up to homotopy), one must show that the image of this hull is exactly the dg-subcategory whose homotopy category is $  \varepsilon_X(\aA)\boxtimes \varepsilon_Y(\bB)  $. Via the Eilenberg–Zilber map for Čech complexes the quasi-equivalence of dg-categories
$$\mu:\perf(\dD_X^{op}\otimes\dD_Y)\xrightarrow{\simeq}\dD_{X\times Y}$$
sends the elementary tensor $  a^{op}\otimes b  $ (with $  a\in\Ob\aA  $, $  b\in\Ob\bB  $) to the external tensor product $  \varepsilon_X(a)\boxtimes\varepsilon_Y(b)  $ in $  D^b(\Coh(X\times Y))  $, hence the image of $  \widehat{\aA^{op}\otimes\bB}  $ under $  \mu  $ is precisely the smallest pretriangulated full sub-dg-category of $  \dD_{X\times Y}  $ whose homotopy category consists of all objects obtained from external tensor products $  \varepsilon_X(a)\boxtimes\varepsilon_Y(b)  $ by taking finite direct sums, cones, and direct summands. But this is exactly the definition of the dg-category whose homotopy category is $  \varepsilon_X(\aA)\boxtimes\varepsilon_Y(\bB)  $.
Hence the composite
$$\widehat{\aA^{op}\otimes\bB}\hookrightarrow\perf(\dD_X^{op}\otimes\dD_Y)\xrightarrow{\mu}\dD_{X\times Y}$$
restricts to a quasi-equivalence onto the dg-subcategory corresponding to $  \varepsilon_X(\aA)\boxtimes\varepsilon_Y(\bB)  $. Composing with $  \varepsilon_{X\times Y}  $ yields the desired equivalence of underlying homotopy categories:
\[\varepsilon_X(\aA)\boxtimes\varepsilon_Y(\bB)\simeq\Ho(\perf(\aA^{op}\otimes\bB)).\]
\begin{flushright}
    $\square$
\end{flushright}

In consequence of the Morita invariance of Hochschild cohomology we get $\Hoc^\star(\aA^{op}\otimes \bB)\cong \Hoc^\star(\widehat{\aA^{op}\otimes \bB})$ and so,  :  

\begin{corollary}In the setting of the previous corollary :  \[\Hoc^\star (\aA^{op}\otimes\bB)\cong \Hoc^\star (\varepsilon_X(\aA)\boxtimes \varepsilon_Y(\bB))\]  
\end{corollary}

Applying our Künneth formula (Corollary \ref{kunnethhochcoh}) we finally conclude to the following result :  

\begin{theorem}   
Let $X$ and $Y$ be smooth and proper projective schemes over $\mathbb{C}$ whose derived category of coherent sheaves respectively contains an admissible quasi-phantom category $A_X$ and $A_y$. Then, the quasi-phantom $A_X\boxtimes A_Y$ of $D^b(X\times Y)$ has Hochschild cohomology given by\[\Hoc^\star(A_X\boxtimes A_Y)\cong \Hoc^{\star}(A_X)\otimes \Hoc^\star(A_Y)\]

In particular, it is non-trivial if $\Hoc^\star(A_X)$ and $\Hoc^\star(A_Y)$ are both non-trivial.
\end{theorem}

\proof We only need to check the hypothesis of theorem \ref{kunnethhochcoh} are fulfilled i.e to show that the underlying algebra of one of our two non-commutative schemes is finite dimensional (in fact both). This is the case because those two are geometric (i.e representing actual smooth and proper projective schemes. The finite dimensionality properties of the algebra of endomorphisms of a generator for $D^b(\Coh(X))$ have been widely studied, notably in \cite{ORLOV2020107096}. 

\begin{flushright}
    $\square$
\end{flushright}

\bibliographystyle{alpha}
\bibliography{biblio}

@article{Gorchinskiy_2013,
   title={Geometric Phantom Categories},
   volume={117},
   ISSN={1618-1913},
   url={http://dx.doi.org/10.1007/s10240-013-0050-5},
   DOI={10.1007/s10240-013-0050-5},
   number={1},
   journal={Publications mathématiques de l’IHÉS},
   publisher={Springer Science and Business Media LLC},
   author={Gorchinskiy, Sergey and Orlov, Dmitri},
   year={2013},
   month=feb, pages={329–349} }

@misc{toen2006homotopy,
      title={The homotopy theory of dg-categories and derived Morita theory}, 
      author={B. Toen},
      year={2006},
      eprint={math/0408337},
      archivePrefix={arXiv},
      primaryClass={math.AG}
}

@article{Kuznetsov_2013,
   title={Height of exceptional collections and Hochschild cohomology of quasiphantom categories},
   volume={2015},
   ISSN={1435-5345},
   url={http://dx.doi.org/10.1515/crelle-2013-0077},
   DOI={10.1515/crelle-2013-0077},
   number={708},
   journal={Journal für die reine und angewandte Mathematik (Crelles Journal)},
   publisher={Walter de Gruyter GmbH},
   author={Kuznetsov, Alexander},
   year={2013},
   month=sep, pages={213–243} }

@article{Orlov_2016,
   title={Smooth and proper noncommutative schemes and gluing of DG categories},
   volume={302},
   ISSN={0001-8708},
   url={http://dx.doi.org/10.1016/j.aim.2016.07.014},
   DOI={10.1016/j.aim.2016.07.014},
   journal={Advances in Mathematics},
   publisher={Elsevier BV},
   author={Orlov, Dmitri},
   year={2016},
   month=oct, pages={59–105} }

@incollection{toen:hal-00772841,
  TITLE = {{Lectures on dg-categories.}},
  AUTHOR = {Toen, Bertrand},
  URL = {https://hal.science/hal-00772841},
  BOOKTITLE = {{Topics in Algebraic and Topological K-Theory}},
  PUBLISHER = {{Springer Berlin Heidelberg}},
  SERIES = {Lecture Notes in Mathematics 2008},
  PAGES = {243-301},
  YEAR = {2011},
  MONTH = Jan,
  DOI = {10.1007/978-3-642-15708-0},
  HAL_ID = {hal-00772841},
  HAL_VERSION = {v1},
}

@misc{toen2007moduli,
      title={Moduli of objects in dg-categories}, 
      author={B. Toen and M. Vaquie},
      year={2007},
      eprint={math/0503269},
      archivePrefix={arXiv},
      primaryClass={math.AG}
}

@misc{Idrissi2022,
      title={Introduction à la théorie de l'homotopie}, 
      author={Najib Idrissi},
      year={2022},
        url={https://cdn.idrissi.eu/main/teaching/homotopie.pdf}}

@misc{tabuada2004une,
      title={Une structure de categorie de modeles de Quillen sur la categorie des dg-categories}, 
      author={Goncalo Tabuada},
      year={2004},
      eprint={math/0407338},
      archivePrefix={arXiv},
      primaryClass={math.KT}
}

@misc{peroux2021coalgebrasdwyerkanlocalizationmodel,
      title={Coalgebras in the Dwyer-Kan localization of a model category}, 
      author={Maximilien Péroux},
      year={2021},
      eprint={2006.09407},
      archivePrefix={arXiv},
      primaryClass={math.AT},
      url={https://arxiv.org/abs/2006.09407}, 
}

@misc{neemansurvey,
      title={The K-theory of Triangulated Categories}, 
      author={Amnon Neeman},
      year={},
      url={https://webhomes.maths.ed.ac.uk/~v1ranick/papers/neemsurv.pdf}, 
}

@misc{shklyarov2007serredualitycompacthomologically,
      title={On Serre duality for compact homologically smooth DG algebras}, 
      author={D. Shklyarov},
      year={2007},
      eprint={math/0702590},
      archivePrefix={arXiv},
      primaryClass={math.RA},
      url={https://arxiv.org/abs/math/0702590}, 
}

@misc{kontsevich2024notesainfinityalgebrasainfinity,
      title={Notes on A-infinity algebras, A-infinity categories and non-commutative geometry. I}, 
      author={Maxim Kontsevich and Yan Soibelman},
      year={2024},
      eprint={math/0606241},
      archivePrefix={arXiv},
      primaryClass={math.RA},
      url={https://arxiv.org/abs/math/0606241}, 
}

@book{hovey2007model,
  title={Model Categories},
  author={Hovey, M.},
  isbn={9780821843611},
  lccn={98345399},
  series={Mathematical surveys and monographs},
  url={https://books.google.fr/books?id=cyn0BwAAQBAJ},
  year={2007},
  publisher={American Mathematical Society}
}

@misc{canonaco2018uniqueness,
      title={Uniqueness of dg enhancements for the derived category of a Grothendieck category}, 
      author={Alberto Canonaco and Paolo Stellari},
      year={2018},
      eprint={1507.05509},
      archivePrefix={arXiv},
      primaryClass={math.AG}
}

@misc{canonaco2022internal,
      title={Internal Homs via extensions of dg functors}, 
      author={Alberto Canonaco and Paolo Stellari},
      year={2022},
      eprint={1312.5619},
      archivePrefix={arXiv},
      primaryClass={math.AG}
}

@misc{bondal2002generatorsrepresentabilityfunctorscommutative,
      title={Generators and representability of functors in commutative and noncommutative geometry}, 
      author={Alexei Bondal and Michel Van den Bergh},
      year={2002},
      eprint={math/0204218},
      archivePrefix={arXiv},
      primaryClass={math.AG},
      url={https://arxiv.org/abs/math/0204218}, 
}

@misc{angel2017bvalgebra,
      title={A BV-algebra Structure on Hochschild Cohomology of the Group Ring of Finitely Generated Abelian Groups}, 
      author={Andrés Angel and Diego Duarte},
      year={2017},
      eprint={1704.03075},
      archivePrefix={arXiv},
      primaryClass={math.RA}
}

@article{LE20141463,
title = {On the Hochschild cohomology ring of tensor products of algebras},
journal = {Journal of Pure and Applied Algebra},
volume = {218},
number = {8},
pages = {1463-1477},
year = {2014},
issn = {0022-4049},
doi = {https://doi.org/10.1016/j.jpaa.2013.11.029},
url = {https://www.sciencedirect.com/science/article/pii/S0022404913002338},
author = {Jue Le and Guodong Zhou}
}

@book{loday93cyclic,
  title={Cyclic Homology},
  author={Loday, J.L.},
  isbn={9783662113899},
  lccn={97022418},
  series={Grundlehren der mathematischen Wissenschaften},
  url={https://books.google.fr/books?id=RzHqCAAAQBAJ},
  year={1993},
  publisher={Springer Berlin Heidelberg}
}

@article{keller:hal-03280332,
  TITLE = {{Hochschild (Co)homology and Derived Categories}},
  AUTHOR = {Keller, Bernhard},
  URL = {https://hal.sorbonne-universite.fr/hal-03280332},
  JOURNAL = {{Bulletin of the Iranian Mathematical Society}},
  PUBLISHER = {{Springer}},
  YEAR = {2021},
  MONTH = Apr,
  DOI = {10.1007/s41980-021-00556-0},
  HAL_ID = {hal-03280332},
  HAL_VERSION = {v1},
}

@article{ORLOV2020107096,
title = {Finite-dimensional differential graded algebras and their geometric realizations},
journal = {Advances in Mathematics},
volume = {366},
pages = {107096},
year = {2020},
issn = {0001-8708},
doi = {https://doi.org/10.1016/j.aim.2020.107096},
url = {https://www.sciencedirect.com/science/article/pii/S0001870820301225},
author = {Dmitri Orlov}
}

@misc{Kontsevitch,
      title={Non commutative motives, talk at the conference on the occasion of the sixty first birthday of Pierre Deligne},
      author={Maxim Kontsevitch},
      year={17-20 oct 2005},
      link={https://www.ihes.fr/~maxim/TEXTS/ncmotives%20(Skoda%20notes).pdf}
}

@article{Orlov_2003,
doi = {10.1070/IM2003v067n03ABEH000437},
url = {https://dx.doi.org/10.1070/IM2003v067n03ABEH000437},
year = {2003},
month = {jun},
publisher = {},
volume = {67},
number = {3},
pages = {535},
author = {D O Orlov},
title = {Quasi-coherent sheaves in commutative and non-commutative geometry},
journal = {Izvestiya: Mathematics}
}

@article{Rosenberg98,
title = {Noncommutative schemes},
author = {Alexander L. Rosenberg},
journal = {Compositio Mathematica},
volume = {112},
pages = {93–125},
year = {1998},
doi = {https://www.cambridge.org/core/journals/compositio-mathematica/article/noncommutative-schemes/A310130B5267707504540607471A5740},
url = {https://www.cambridge.org/core/services/aop-cambridge-core/content/view/A310130B5267707504540607471A5740/S0010437X98000372a.pdf/noncommutative_schemes.pdf},
}

@misc{keller2006differentialgradedcategories,
      title={On differential graded categories}, 
      author={Bernhard Keller},
      year={2006},
      eprint={math/0601185},
      archivePrefix={arXiv},
      primaryClass={math.KT},
      url={https://arxiv.org/abs/math/0601185}, 
}

@article{Keller1994,
     author = {Keller, Bernhard},
     title = {Deriving {DG} categories},
     journal = {Annales scientifiques de l'\'Ecole Normale Sup\'erieure},
     pages = {63--102},
     publisher = {Elsevier},
     volume = {Ser. 4, 27},
     number = {1},
     year = {1994},
     doi = {10.24033/asens.1689},
     mrnumber = {95e:18010},
     zbl = {0799.18007},
     language = {en},
     url = {http://www.numdam.org/articles/10.24033/asens.1689/}
}

@misc{kuznetsov2009,
      title={Hochschild homology and semiorthogonal decompositions}, 
      author={Alexander Kuznetsov},
      year={2009},
      eprint={0904.4330},
      archivePrefix={arXiv},
      primaryClass={math.AG},
      url={https://arxiv.org/abs/0904.4330}, 
}

@misc{orlov2008remarksgeneratorsdimensionstriangulated,
      title={Remarks on generators and dimensions of triangulated categories}, 
      author={Dmitri Orlov},
      year={2008},
      eprint={0804.1163},
      archivePrefix={arXiv},
      primaryClass={math.AG},
      url={https://arxiv.org/abs/0804.1163}, 
}

\end{document}